\documentclass[3p]{elsarticle}
\usepackage[utf8]{inputenc}
\usepackage[]{todonotes}
\usepackage[ruled]{algorithm2e}
\usepackage{endnotes}
\usepackage{pgfplots}
\usepackage{lscape}
\usepackage{setspace}
\pgfplotsset{compat=newest}

\usepackage{rotating}
\usepackage{float}
\usepackage{verbatim}

\usepackage{lscape}
\usepackage{amsmath,amsfonts,amsthm,amssymb}
\def\setR{\mathbb{R}}

\usepackage{footnote}
\makesavenoteenv{tabular}

\usepackage{hyperref}
\usepackage{comment}
\usepackage{graphics}
\usepackage{graphicx}
\usepackage{adjustbox}
\usepackage{tikz}
\usepgfplotslibrary{dateplot}
\usepgfplotslibrary{groupplots}
\usepackage{caption}
\usepackage{subcaption}
\usepackage{psfrag}
\usepackage{soul}
\usepackage{mwe}
\usepgfplotslibrary{fillbetween}
\usetikzlibrary{er,positioning, calc, shapes, arrows, chains, fit}
\definecolor{myblue}{RGB}{19,145,215}
\definecolor{mygreen}{RGB}{80,176,50}
\definecolor{myred}{RGB}{222,36,16}
\usepackage{fancyhdr}
\usepackage{rotating}

\begin{document}

\begin{frontmatter}
\title{An actor-critic algorithm with policy gradients to solve the job shop scheduling problem using deep double recurrent agents} 

\author[1]{Marta Monaci}
\ead{marta.monaci@uniroma1.it}

\author[1]{Valerio Agasucci}
\ead{agasucci@diag.uniroma1.it}

\author[2]{Giorgio Grani$^{*,}$\footnote{$^*$ Corresponding author. }}
\ead{g.grani@uniroma1.it}

\address[1]{Sapienza University of Rome, Dep. of Computer Science, Control and Management Engineering, Rome, Italy
}

\address[2]{Sapienza University of Rome, Dep. of Statistics, Rome, Italy
}

\begin{abstract}

In this work, we present a method that applies Deep Reinforcement Learning, an approximate dynamic programming procedure using deep neural networks, to the job shop scheduling problem (JSSP).  The aim is to show that a greedy-like heuristic trained on a subset of problems, can effectively generalize to some extent to unseen instances, and be competitive compared to other methods.

We model the JSSP as a Markov Decision Process and we exploit the efficacy of reinforcement learning to solve the problem.  { We adopt an actor-critic scheme based on policy gradients, specifically the Proximal Policy Gradient method, where the action taken by the agent is influenced by policy considerations on the state-value function.} The procedures take into account the challenging nature of JSSP, where the state and the action space change for every instance and after each decision. To tackle this variability, we introduced a novel model based on two incident Long-Short Term Memory networks, followed by an encoding model, different in structure for both the actor and the critic.
 
Experiments show the algorithm reaches good solutions in a short time,  proving that is possible to generate new greedy heuristics just from learning-based methodologies.   We compared our algorithms against several established heuristics, an adaptive method, a commercial solver based on branch and cut, and another approach based on Deep Reinforcement Learning, proving the validity of the proposed method in terms of time and makespan.
The model can generalize, to some extent, to larger problems originating from a different distribution.
\end{abstract}

\end{frontmatter}

\parindent 0cm

\pagestyle{fancy}
\lhead{}
\renewcommand{\headrulewidth}{0pt}
\rhead{\textcolor[gray]{0.3}{Published at: \href{https://doi.org/10.1016/j.ejor.2023.07.037}{https://doi.org/10.1016/j.ejor.2023.07.037}}}
\thispagestyle{fancy}  

\textcopyright{ 2023. Licensed under the Creative Commons  \href{https://creativecommons.org/licenses/by/4.0/}{CC-BY 4.0}.}

\section{ Introduction}
\label{Introduction}

There is a large amount of research being done in the fields of Operations Research and Mathematical optimization to solve problems   {related} to scheduling. The ability to schedule tasks efficiently permeates a lot of human activities, from classical industrial or transportation applications to services, computer science,  and even to healthcare.  The problem is so   {common that} we encounter its solutions on a daily basis, like when consulting the board to choose which train to take, or on a website deciding a flight for our next vacation, or    {when ordering  online} and it gets assembled together in one order and delivered to us. Basically, every processed good on the market is a result of a set of operations,   {scheduled to} get the final product. Every respectful scheduling problem is characterized by a set of procedures, called tasks or operations, that has to be executed in order to complete a job. Depending on the structure of this set (e.g. ordered or not) and where the task has to be exploited, for instance on a fixed location or dislocated points, we have a different type of scheduling problem. Other than that, one may add conditions related to the release or due date of an order, the time window for which a certain processing point is available, set-up times, or every other kind of business requirement. In this paper, we will address the well-known minimum makespan Job Shop Scheduling problem (JSSP), where we have a certain number of jobs needing to be completed, and each job is composed of a list of operations that can be exploited in specific processing points, also referred as machines. {The aim is to schedule the operations on the machines in order to minimize the makespan, which is the delay accumulated by the entire system}. 

\medskip
In 1977, Lenstra et al. showed in \cite{Lenstra1977Complexity} that the JSSP is NP-hard and, as reported in \cite{Blazewicz1996Thejob}, only a few specific cases are polynomially solvable. 
For instance, job shop problems with two jobs are efficiently solved by the geometric approach, described in \cite{Brucker1988Anefficient} and originally presented in \cite{Akers1956Letter}. In \cite{Johnson1953Optimaltwo}, Johnson presented a simple decision rule to optimally solve the two-machine flow shop problem, in which each job must be executed on the machines in the same order. 
In order to extend this result, Jackson et al. \cite{Jackson1959Anextension} proposed an efficient procedure to solve the two-machine job-shop problem where each job is composed of at most two tasks by reducing it to the flow shop case (see also \cite{Sotskov1991Thecomplexity}). 
Then, in \cite{Hefetz1982AnEfficient}, an efficient algorithm to solve two-machine job-shop problems with unit processing times is presented. The number of steps and the storage space of the method is linear in the total number of operations, thus proving that the problem belongs to the \textit{P}-class. 
Later, Brucker \cite{Brucker1994Apolynomial} developed a polynomial algorithm to solve the two-machine job shop case with a fixed number of jobs, even if machine repetition is allowed. In \cite{brucker2007complexity}, the NP-hardness of scheduling problems with a fixed number of jobs was investigated.
Slight modifications to these categories of problems have been shown to be hard to solve. For instance, three-machine job shop problems where each job has at most two operations or with a fixed number of three jobs are NP-hard (\cite{Gonzalez1978Flowshop}, \cite{Sotskov1995NP}), (\cite{lenstra1979computational}). In \cite{Garey1976complexity}, it is shown that job shop scheduling problems with more than two machines are NP-complete.





Therefore, solving a JSSP instance can be, in general, quite complicated. Lots of researchers tried their best to build nice and performing algorithms, investigating tons of different optimization strategies. It would be impossible for us to cover even only the most important papers in this context, so we will just take a glimpse of the overall picture, redirecting to \cite{zhang2019review} and \cite{chaudhry2016research} for a couple of comprehensive surveys. In principle, we divide the class of algorithms between exact and heuristic methods. The first group concentrates on finding the optimal solution, the best schedule above all, according to the objective. They should be used when the size of the problem and the time required to find a solution are compatible with the application. The second group focuses on finding something good enough in a reasonable short time,  respecting feasibility, but with no guarantee of finding the optimum. 
Exact algorithms are largely associated with mathematical integer programming. In literature a multitude of approaches has been applied, like column-generation in \cite{gelinas2005dantzig} and \cite{lancia2007compact}, or branch-and-bound in \cite{Carlier1989Analgorithm} and \cite{Brucker1994Abranch}, branch-and-cut in \cite{karimi2015lot}, and Lagrangian relaxations as in \cite{hoitomt1990lagrangian} and \cite{chen2003alternative}. On the heuristic side of the palisade, we can go through all kinds of known methods: simulated annealing, genetic algorithms, local search, tabu search, and so on. A nice survey on swarm and genetic algorithms can be found in \cite{gao2019review}, while one on meta-heuristics in \cite{mhasawade2017survey}. Also, dynamic programming plays a big role, but the complexity of the instances makes it useful only in an approximate way. A review of dynamic programming and JSSP is presented in \cite{mohan2019review}.
\medskip

With more data being generated and more and more operations being automated, machine learning-based techniques have found their place within the JSSP community. This trend is common to other aspects of optimization,  as described in \cite{bertsimas2019machine} and \cite{bengio2020machine} in the case of machine learning, deep learning and combinatorial optimization. Every time it is common to encounter some kind of approximation or when you need to solve several instances of the same type sequentially, machine learning may play a role.  The possible ways of the interconnection of the two fields are countless, with researchers bringing new ideas day by day. Sometimes, it is useful to approximate reality or complex systems by using surrogate functions and inserting them directly into the formulation. The concept of a surrogate model has been used for a long time in optimization, in the form of linear or quadratic approximations, and it is now extended to more complex models. An example of using neural networks to ease the, otherwise unbearable, complexity of the involved physics in manufacturing can be found in  \cite{pfrommer2018optimisation}, whereas in \cite{hottung2020learning} there is a usage of autoencoders to detect useful embeddings for routing problems.
On other occasions, we may use machine learning inside a complex algorithmic framework, improving some of the frequent decisions taken through the process. For instance, we may estimate a score suggesting the next variable to choose for branching, or which of the several local heuristics fits more in the current iteration, or the most promising cuts to add. As anticipated, a great deal has been achieved in the mixed-integer programming domain. Just to mention a few works, in \cite{khalil2016learning} the authors make use of support vector machines to mimic strong branching, whereas in \cite{gupta2020hybrid} they exploit a similar task but using deep neural networks. Authors in  \cite{tang2020reinforcement} chose the most promising Gomory's cuts (see \cite{wolsey1999integer}) by means of long-short term memory networks and reinforcement learning. Recently, in \cite{nair2020solving}, scientists from DeepMind and Google Research proposed a method combining deep branching and deep diving in a branch-and-cut algorithm, taking advantage of supervised learning, graph convolutional neural networks, and alternating descent method of multipliers.
Another important aspect is to learn primal heuristics that return directly to near-optimal solutions.  In \cite{bertsimas2019online} and \cite{bertsimas2021voice}, this is done opportunely by studying some properties of the optimization problem (the so-called \emph{voice of optimization}), then they train a neural network to directly produce a solution, and finally, they use a projection technique to make it feasible, leading to an extremely fast heuristic. For stochastic problems, in \cite{bengio2020learning} the authors did something similar, whereas in \cite{cauligi2020learning} this methodology has been used to generate robust trajectories in planning robot agents. Finally, in \cite{hottung2020deep}, this concept is extended to approximate bounds and inserted within a tree search for solving the container pre-marshaling problem.

In some cases, one may learn a new algorithm from scratch thanks to reinforcement learning (RL). We will also adopt RL through this paper, and it will be described in detail later on. For now, consider it to be a greedy-like method, where decisions are taken at each step by a learning-based operator. Just to cite a few papers, in \cite{khalil2017learning},  \cite{gasse2019exact} and \cite{drori2020learning} the use of RL to solve directly optimization problem is exploited with remarkable results. The three works differ in the type of RL framework adopted (Q-learning, policy gradient, or actor-critic) and in the neural network structure.  On the same page, in \cite{agasucci2020solving}, we used Deep Q-learning to solve the train dispatching problem, comparing two approaches: a centralized one looking at the overall rail network, and another decomposing the problem by train with a limited view of the surroundings.
\medskip

After this brief description of machine learning and optimization, we are ready to go back  and discuss the intersection with JSSP. The problem has received interest from the machine learning community for a long time, as supported in \cite{ccalics2015research}. This was before the deep learning era, so before efficient fast-computing libraries were established as a standard for neural networks.

No more than ten years ago, scientists and professionals started realizing the great potential of deep learning to solve complex tasks. When we talk about deep learning, we refer to everything connected to neural networks with more than one hidden layer. In \cite{goodfellow2016deep}, one may find a larger description of the concepts and methodologies of deep learning.

Most of the approaches linked to JSSP are also connected to RL or to approximate dynamic programming in general.  A remarkable example can be found in \cite{Sotskov1996Adaptive}, where the authors studied the strength of adaptive algorithms, taking advantage of a particular graph structure to ease the computational effort. The method was later extended to the case of parallel machines JSSP in \cite{gholami2014solving}. This methodology shares a lot of elements with the RL approaches, including training. In section \ref{sec:computational_experience}, we show a comparison between our proposed method and the adaptive algorithm, showing how RL approaches with Deep neural networks are able to outperform the other on the proposed instances. Another interesting work can be found in \cite{zhang1995reinforcement}, where the authors present  {an innovative but shallow RL framework}, to solve an extension of the JSSP for a real application at NASA. {In particular, they proposed a RL algorithm, based on Temporal Difference Learning, aimed at learning ad-hoc repair heuristics to produce good conflict-free schedules.}    { In \cite{zhang2020learning}, the authors used a Deep RL apporach based on Graph Neural Networks and Proximal Policy Optimization. We compared ourselves with the work presented in \cite{zhang2020learning}.  In \cite{tassel2022reinforcement} the authors used Neural Networks in combination with a RL method based on Natural Evolution Strategies \cite{wierstra2014natural}, to solve dispatching problems on an industrial scale. }

\medskip

This paper describes our proposal to tackle the JSSP using RL. The algorithm adapts the Policy Proximal Optimization (PPO) algorithm firstly presented in \cite{schulman2017proximal} to the JSSP, making use of a suitable representation of the environment as a Markov Decision Process (MDP). PPO belongs to the family of actor-critic RL algorithms, for which we developed two special deep neural networks for both the actor and the critic, based on two concatenated Long Short-Term Memory networks (LSTMs). This architecture is proven to be effective and flexible to the number of jobs, operations, and machines. We compared our algorithm against standard JSSP heuristics (as resumed in \cite{Panwalkar1977ASurvey}), the adaptive algorithm presented in \cite{Sotskov1996Adaptive}, the branch-and-cut algorithm implemented in the known solver CPLEX from IBM (\cite{cplex}),    { and the Deep Reinforcement learning algorithm described in \cite{zhang2020learning}}, obtaining good average results.\medskip

The paper is organized as follows: in Section \ref{sec:preliminaries_and_notation} we formalize the JSSP and we give the basic elements and ideas behind the paradigm of reinforcement learning, with a special focus on actor-critic methods. In Section \ref{sec:jssp_as_mdp}, we describe the JSSP as a Markov decision process,  making the problem solvable by using RL. Then, in Section \ref{sec:agents_deep_models}, we introduce the deep neural networks used as a learning model in our algorithm. Finally, in Section \ref{sec:computational_experience}, we illustrate the experiment.

\section*{Our contribution}

This work presents a policy proximal optimization algorithm with deep agents to tackle the JSSP. The major findings in this paper are:
\begin{itemize}

\item A \textbf{novel model} both for the actor and the critic, using two concatenated Long Short-Term Memory networks (LSTMs).
 \item The method is \textbf{flexible} and not related to a single application. In particular, the Double LSTM structure allows to vary arbitrarily the number of jobs, operations and machines adopted.
 \item The models \textbf{generalize} to some extent for larger and more complex instances maintaining good solution quality.

\item The computational experience is conducted both on time and solution quality \textbf{against a commercial solver, an adaptive heuristic, 17 rule-based heuristics, and a    {Deep Reinforcement Learning approach}}, showing good average results.
   
\end{itemize}

\section{Preliminaries and notation}\label{sec:preliminaries_and_notation}

In this section, we describe the Job Shop Scheduling problem (JSSP) as an optimization problem and we enter the world of reinforcement learning (RL), actor-critic methods and the Proximal Policy optimization algorithm (\cite{schulman2017proximal}).

\subsection{The Job Shop Scheduling Problem}
Scheduling is a decision-making process finding a temporal allocation of shared and limited resources to activities to optimize some desired objective. In this project, we are tackling the {n$\times$m minimum makespan} JSSP, denoted by $Jm\ ||\ C_{max}$ according to the three-field notation introduced by Graham et al. \cite{GrahamOptimization1979}. It will be referred to as JSSP throughout the rest of the paper without loss of information.

In its standard form, we have a bunch of workers and a set of operating stations doing some service. Each worker has a good to process, and a list of ordered required operations to be performed to get the final product. The point of JSSP is determining the exact timing for which each worker should go to an operative station and perform some kind of process on the good, minimizing the overall time for all the workers. From now on the operating stations will be called \textit{machines}, the processes \textit{tasks} (or operations) and the list of ordered operations \textit{jobs}.

Let $\mathcal{J}=\{j\}_{j=1}^n$ be the set of {jobs}, which has to be processed on the set $\mathcal{M}=\{k\}_{k=1}^{m}$ of {machines}.
Each job $j$ has a given processing sequence of $n_j$ different machines, with $n_j \leq m$. 
A task (or operation) is the activity that job $j\in\mathcal{J}$  must execute on  machine $k\in\mathcal{M}$, and it is denoted by the pair $(j,k)$. 

Therefore, each job is  a list of different tasks to be performed. We adopt the notation for which $(j,k) \prec (j,h),$
means that operation $(j,k)$, the $i-th$ operation of job $j$, precedes operation $(j,h)$, the $(i+1)$-th operation of job $j$, $\forall i = 1,\dots,n_j-1$. These rules define the {precedence constraints} for the problem. A processing time $p_{jk}$ is associated to each operation $(j,k)$. The set of all the operations is denoted by $\mathcal{O}$.

We assume to be in a {no-preemption} regime so that operations can not be interrupted. Moreover, each machine can not process more than one job at the same time, meaning there is no overlap. We call  $t_{jk}$ the starting time of the operation $(j,k)$, and  its
 {completion time} $C_{jk}$, which is the time interval elapsing from the start of the whole process to the execution of the operation itself, i.e. 
$C_{jk}=t_{jk}+p_{jk}$.
The optimum is reached by minimizing the {makespan}, denoted by $C_{max}$, which is the maximum completion time of all the operations, i.e. $ C_{max}=\max_{(j,k)\in \mathcal{O}} C_{jk}$.

\subsection{Reinforcement Learning}

Reinforcement learning (RL) is a paradigm of machine learning, alongside supervised and unsupervised learning. There are four elements in an RL framework: agent, action, state, and reward. They all operate imitating the decision process in a real-world setting. The agent is the decision-maker, the actions are the set of options it is allowed to do, the state is an encoding of the environment it operates into, and the reward is what it gets after making an action. At every step, the agent observes the state, takes an action, and waits for the environment to return to its new form, the next state, alongside a reward for selecting that action. After that, the system is ready for  a new iteration. In RL, the agent learns from its own choices, step by step, game by game, self-generating data. The goal of the agent is specified by an objective function, dependent on the collected rewards.

RL can be formalized as a Markov Decision process, see \cite{sutton2018reinforcement}, and it can be seen as an approximate dynamic programming method, see \cite{bertsekas1995dynamic} and \cite{bertsekas2019reinforcement}).
Following the dynamic programming terminology, the goal of RL is to learn an optimal strategy, called {policy}, allowing the agent to solve the problem by maximizing its total reward. Unfortunately, we are not able to inspect the full tree of possibilities and alternatives an agent may encounter, since the size of the decision tree would be, for NP-hard problems like JSSP, unbearable for standard computational resources. For this reason, we discuss policies maximizing the expected total reward, so that we infer optimality by just observing a portion of the decision space.
\medskip

An episode is an instance to be solved by an RL algorithm while training the agent. To complete an episode, the RL method goes step by step until reaching the $T$-th state, where $T$ is the last possible iteration. Steps within an episode are indexed by $t=0,1,\dots,T$, where $t=0$ is the initial state. We use the notation $s_t$ for states, $a_t$ for actions and $r_t$ for rewards, all depending on $t$. The set of all the possible states is $\mathcal{S}$,  whereas 
 $\cal A$ is the set of all the possible actions and $ R: \mathcal{S}\times \mathcal{A} \mapsto \setR$ is a function associating states and actions to rewards.
 A  {stochastic} policy is identified by  the  function $\pi(a_t|s_t)$, representing the probability to take the action $a_t$ when the state $s_t$ is observed where {$\lim_{t \rightarrow \inf} \pi(a_t,s_t)$ will tend to be an optimal, deterministic policy.}

The expected cumulative reward is the objective function to be maximized in the RL framework.
\begin{equation}\label{eq:cumulative_reward}
    J = \mathbb{E}_{\pi}\left[\sum_{t=1}^T r_t\right],
\end{equation}
where $\mathbb{E}_\pi$ is the expected value computed according to the policy distribution $\pi$.

Several RL methods make use of the state-value function, which is a measure of the expected cumulative reward from a certain step $k$ (also called expected reward-to-go) when observing a state $s_k$.
$$
     V_\pi(s_k) = \mathbb{E}_\pi\left[\sum_{t= k+1}^T  r_t \Big\vert \mathcal{S}_k = s_k \right]
$$

For similar purposes, it is handy to define the advantage function, which evaluates the expected improvement when selecting an action,taking as input an action $a_k$ and the state $s_k$.  This can be thought of as a sort of differential measure over the reward-to-go, see \cite{PaperGAE}.  We present a version of the formula valid in the context of JSSP, where it is always possible to derive the reward $r_k$ from the state and the taken action.
\begin{equation}\label{eq:advantage}
     A_\pi(a_k, s_k)  = r_{k} + V_\pi (s_{k+1}) - V_\pi(s_k)
\end{equation}

Among the numerous classes of RL algorithms available in the literature, we are going to focus on actor-critic methods, which are characterized by having the agent separated into two decision entities: the actor and the critic. The reason behind this separation is to allow a policy improvement through an estimation of the state-value function, combining both value-based and policy-improvement algorithms. The critic approximates the state-value function  $\hat{V}(s)$, while the actor updates and improves a model of the stochastic policy $\hat{\pi}$  by taking into account the critic estimation while maximizing the total expected reward.

Several algorithms in the actor-critic sense have been proposed over the last years, and the main difference lies in how to properly train the models describing the actor and the critic. For instance, the update rule may be based on the Bellman equation (see \cite{bellman1966dynamic}), as proposed in \cite{barto1983neuronlike}, \cite{konda1999actor} and more recently in \cite{lawhead2019bounded}. In this paper, we will focus on actor-critic methods basing their updating operations on policy gradient rules. One of the first examples can be found in \cite{konda2003onactor}, where a class of two time-scale algorithms is presented, in which the critic uses classical temporal difference learning, and the actor uses policy gradient based on the critic estimation. More recently, in \cite{kakade2001natural} and \cite{kakade2002approximately}, the authors presented some properties for a specific class of policies, allowing for the definition of a lower bound over the difference between two policies. The update is obtained by minimizing this bound so that the new policy will improve as far as known using the previous information. This idea is extended in the Trust Region Policy Gradient (TRPO) \cite{schulman2015trust} to more general policies, and an approach based on a Kullback Liebler divergence trust region is shown to work efficiently. Finally, in \cite{schulman2015trust}, the trust region approach is abandoned to an unconstrained one in Proximal Policy Gradient (PPO). Since in the following of the paper we will discuss the PPO, we will now   {discuss in more detail}    {proximal policy gradient based actor-critic algorithms.}

In the following, we will assume that both the actor and the critic models have as weights  $\theta$ and $\omega$, respectively.

For ease of explanation, we will keep this dependence clear when needed in our formulas. For instance, the total expected reward presented in \eqref{eq:cumulative_reward} can be rewritten as 
     $J_{\theta}= \mathbb{E}_{\pi_{\theta}}\left[\sum_{t=1}^Tr_t\right]$, stressing the dependency of the actor weights in the objective function.

To optimize the total reward function \eqref{eq:cumulative_reward}, it is possible to derive an estimator for the gradient, and use it in gradient-based optimization methods. For computational efficiency reasons, we report the formula in terms of log probabilities.

\begin{equation}\label{eq:policy_gradient}
    \nabla_\theta J_\theta \propto  \mathbb{E}_{\pi_\theta} \left( \sum_{a_t} \nabla\log \pi_\theta(a_t|s_t) {A_\omega}(s_t,a_t) \right)
\end{equation}

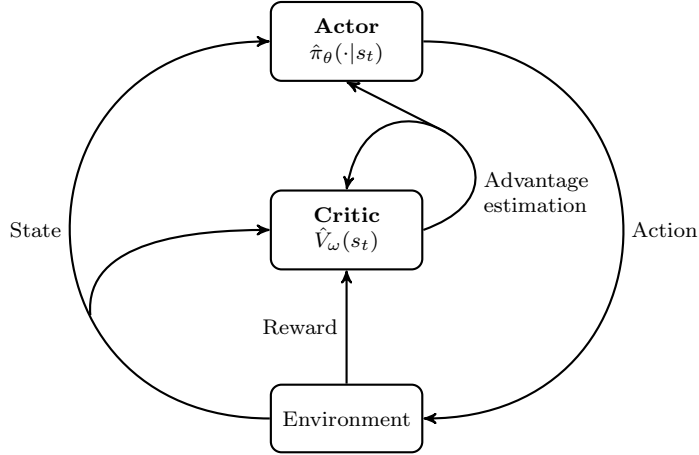
\begin{figure}
    \centering
    \begin{tikzpicture}[->,>=stealth',auto,node distance=1.5cm,
  thick,main node/.style={rectangle, draw,  
    text width=6em, text centered, rounded corners, minimum width=4em, minimum height=3.5em}]
    \footnotesize 
\node[main node] (1) {\textbf{Actor} $\hat{\pi}_\theta(\cdot|s_t)$};
\node[main node] (2) [below of = 1, xshift=-0.0cm, node distance = 2.5cm] {\textbf{Critic} $\hat{V}_\omega(s_t)$};
\node[main node, minimum width=4em,  minimum height=3em] (3) [below of = 2, node distance = 2.5 cm] {Environment} ;

\node (4) [left of = 2, node distance = 3.5 cm] {};


\node (10) [text width=6em, right of = 2, xshift = -0.3cm, yshift = 0.5 cm, node distance = 3 cm] {Advantage estimation};


\draw [->] (1.east) to [out=360,in=360, distance=3.5cm] node {Action} (3.east);

\draw [->] (-3.37,-3.63) to [out=95,in=180, distance=1cm] node {} (2.west);

\draw [->] (3.west) to [out=360,in=360, distance=-3.5cm] node {State} (1.west);

\draw [->] (3.north) to  node {Reward} (2.south);

\draw [-> ] (2.east) to [out=20,in=80, distance=2cm] node {} (2.north);

\draw [-> ] (1.3, -1.2) to  node {} (1.south);



\end{tikzpicture}
    \caption{Actor-critic framework.}
    \label{fig:Actor_Critic}
\end{figure}

The gradient depends on both the actor, $\pi$, and the critic, $A$, so that the log-probability associated with a state-action pair $(s_t,a_t)$ is proportional to the advantage of this pair and thus, the gradient indicates the direction of greatest improvement (locally). Actor-critic algorithms differ in how they compute the estimators, on the version of the policy gradient they adopt, and finally, on the optimization algorithm to improve the expected total reward. The basic scheme of an actor-critic algorithm is illustrated in Figure \ref{fig:Actor_Critic}.
\bigskip

In the PPO, the policy improvement is pursued through a stochastic gradient ascent step over a surrogate function approximating the total expected reward \eqref{eq:cumulative_reward}. 
The loss, called KL-penalized objective, is composed of a surrogate advantage and a Kullback-Leibler divergence multiplied by a negative penalty term $-\beta$. The actor is trained to maximize this objective, refining the approximation from a step to the other, as reported below:

\begin{equation}\label{eq:actor_update}
    \max_{\theta}\ \  {\mathbb{E}}_{\pi_\theta}\left[
     \frac{\pi_\theta(a_t|s_t)}{ \pi_{\theta_{old}}(a_t|s_t)} {A}_{\omega}(s_t,a_t)- \beta D_{KL}\left({\pi}_{\theta} ( \cdot |s_t ),{\pi}_{{\theta}_{old}} ( \cdot | s_t )\right) \right],
\end{equation}

 where $\pi_{\theta_{old}}$ indicates the policy parameters at the previous step, and the parameter $\beta$ influences how much the new policy may diverge from the old one according to the following rule:

\begin{equation}\label{eq:beta_update}
\beta = \left\{
\begin{array}{ll}
    2 \beta_{old}, &\displaystyle \text{if } {\mathbb{E}}_{\pi_{\theta}}\left[D_{KL}\left({\pi}_{{\theta}},{\pi}_{{\theta}_{old}}\right)\right]>1.5\delta \\
   \beta_{old} / 2 ,  &\displaystyle \text{if } {\mathbb{E}}_{\pi_{\theta}}\left[D_{KL}\left({\pi}_{{\theta}},{\pi}_{{\theta}_{old}}\right)\right]<1.5/\delta \\
    \beta_{old},  & \text{otherwise}
\end{array}
\right.,
\end{equation}

where $\delta$ is a target value chosen heuristically. During the algorithm, the penalty coefficient $\beta$ adapts rapidly, and, according to \cite{schulman2017proximal}, its starting value does not affect significantly the training. 

The critic is updated trough the minimization of a mean square error loss function, using  the data collected by the actor in several roll-outs. That is
\begin{equation}\label{eq:critic_update}
    \min_{\omega} {\mathbb{E}}_t\left[\left(V_{\omega}(s_t)-{V_t}^{target}\right)^2\right] =  \min_{\omega} {\mathbb{E}}_t   \left[\left(V_{\omega}(s_t)-R_t\right)^2\right], 
\end{equation}
where the rewards-to-go $R_t$ is the sum of the rewards collected from $t$ to $T$.
\medskip

The algorithm works in the following way: at each episode $k$, it performs $N$ roll-outs, running the policy for $T$ steps,  therefore generating $N\cdot T$ samples.
For each episode $k$,  the rewards-to-go  $R_t$, $t=1,\dots,T$, are computed and stored. 
Then,  the policy is updated solving \eqref{eq:actor_update} with  ADAM, \cite{kingma2014adam}, using the $N \cdot T$ samples. 
Finally, the state-value function is updated according to \eqref{eq:critic_update}, performing some iterations of  ADAM. Algorithm \ref{alg:PPO} reports the procedure.
\medskip

\begin{algorithm}[H]\caption{Proximal Policy Optimization (PPO) with adaptive Kullback-Leibler penalty}\label{alg:PPO}
\footnotesize
	\DontPrintSemicolon
	{\bf Input:} number of roll-outs per episode $N$, 
	termination step $T$, 
	stochastic gradient ascent iterations $L_{\text{actor}}$,
	ADAM iterations $L_{\text{actor}}$,
	target KL divergence $\delta$,
	$\theta = \theta_0$, 
	$\omega = \omega_0$, 
	$\beta = \beta_0$.\;
	\BlankLine
	 
	\For{episode $k = 0,1, \dots$}{
	    \For{roll-out $i = 1,\dots,N$}{
	        Run policy $\pi_{\theta}$ in the environment for $T$ time-steps.\;
        	Compute rewards-to-go $R_1,\dots,R_T$ associated to the $i-th$ roll-out and store them.\;
            Compute advantages associated to the $i-th$ roll-out and store them.\;}
        	Update the actor parameters $\theta$ according to \eqref{eq:actor_update}, performing $L_{\text{actor}}$ iterations of ADAM.\;

    	Update $\beta$ according to \eqref{eq:beta_update}.\;
	Update the critic parameters $\omega$ according to \eqref{eq:critic_update}, performing $L_{\text{critic}}$ iterations of ADAM.\;
	\;}
\end{algorithm}

\section{JSSP as a Markov Decision process}\label{sec:jssp_as_mdp}

JSSP can be solved in several ways, depending on which element we are looking for to find a solution. Our idea is to choose one operation at a time,  deciding whenever a machine is available. In other words, we have a set of jobs, each one being a list of operations to be addressed their specified machine, and we want to find a, possibly good, solution in a greedy fashion. Whenever a machine   {is available to accept jobs}, we select one (and only one) job that can proceed in the queue. In this way, the total number of decisions is equal to the total number of tasks. 

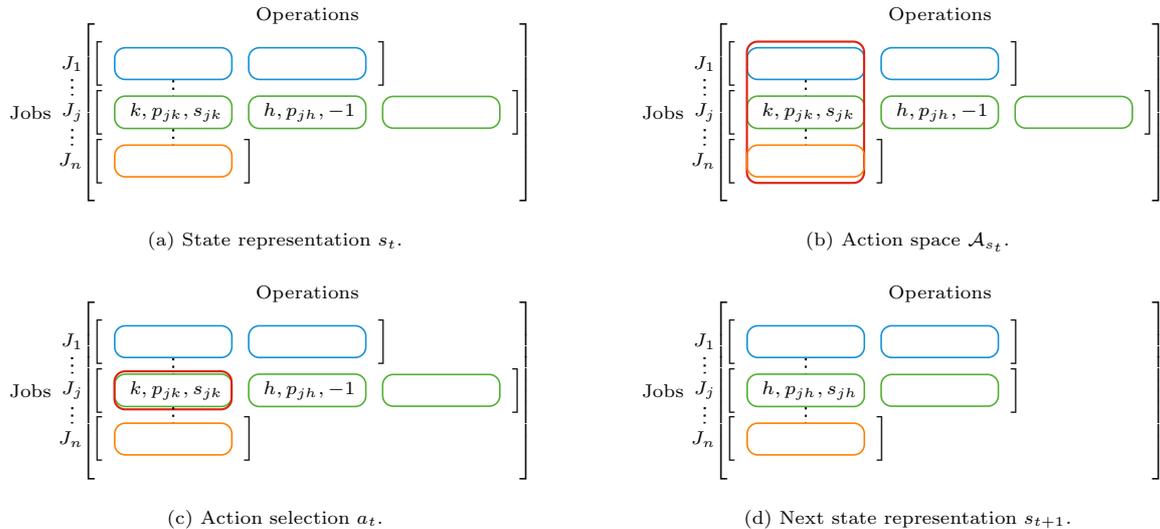
\begin{figure}
\centering 
\scriptsize 

\begin{subfigure}{0.48\textwidth}
\centering 
\begin{tikzpicture}[->,>=stealth',auto,node distance=1.42cm,
  semithick,
  main node blue/.style={rectangle, draw=myblue, text width=4.1em, minimum width = 5.5em, text centered, rounded corners, minimum height=1.5em},
  main node green/.style={rectangle, draw=mygreen, text width=4.1em, text centered,  minimum width = 5.5em,, rounded corners, minimum height=1.5em},
  main node orange/.style={rectangle, draw=orange, text width=4.1em, minimum width = 5.5em, text centered, rounded corners, minimum height=1.5em}
  ]

    \node (generator) {Operations};
\node[main node blue] (4) [below=0.2cm of generator] {};
\node[main node blue] (1) [left=0.2cm of 4] {};
\node[main node green] (2) [below=0.2cm of 1] {$k, p_{jk}, s_{jk}$};
\node[main node orange] (3) [below=0.2cm of 2] {};
\node[main node green] (5) [right=0.2cm of 2] {$h, p_{jh}, -1$};
\node[main node green] (6) [right=0.2cm of 5] {};


\node (left-paren) [left=0.1cm of 2] {$\left[\vphantom{\begin{array}{c}     1\\ 1\\ 1\\ 1\\ 1\\ 1\\ 1\\ 1\\ 1 \end{array}}\right.$};
\node (right-paren) [right=0.1cm of 6] {$\left.\vphantom{\begin{array}{c}   1\\ 1\\ 1\\ 1\\ 1\\ 1\\ 1\\ 1\\ 1 \end{array}}\right]$};

\node (left-paren1) [left=0cm of 2] {$\left[\vphantom{\begin{array}{c}     1\\ 1 \end{array}}\right.$};
\node (right-paren1) [right=0cm of 6] {$\left.\vphantom{\begin{array}{c}     1\\ 1 \end{array}}\right]$};
\node (left-paren2) [left=0cm of 1] {$\left[\vphantom{\begin{array}{c}     1\\ 1 \end{array}}\right.$};
\node (right-paren2) [right=0cm of 4] {$\left.\vphantom{\begin{array}{c}     1\\ 1 \end{array}}\right]$};
\node (left-paren3) [left=0cm of 3] {$\left[\vphantom{\begin{array}{c}     1\\ 1 \end{array}}\right.$};
\node (right-paren3) [right=0cm of 3] {$\left.\vphantom{\begin{array}{c}     1\\ 1 \end{array}}\right]$};


 \node (1b) [left=0.3cm of 1] {$J_1$};
\node (2b) [left=0.3cm of 2] {$J_j$};
\node (3b) [left=0.3cm of 3] {$J_n$};
 \node (00) [left=0.72cm of 2] {Jobs};
\draw [-, dotted, thick] (1.south) to  (2.north);
\draw [-, dotted, thick] (2.south) to  (3.north);
\draw [-, dotted, thick] (1b.south) to  (2b.north);
\draw [-, dotted, thick] (2b.south) to  ($(3b.north)+(0.18mm,0)$);

\end{tikzpicture}
    \caption{State representation $s_t$.}
  \label{fig:state}
\end{subfigure}%
\hfill
\begin{subfigure}{0.48\textwidth}
\centering 
  \begin{tikzpicture}[->,>=stealth',auto,node distance=1.42cm,
  semithick,
  main node blue/.style={rectangle, draw=myblue,text width=4.1em,  minimum width = 5.5em, text centered, rounded corners, minimum height=1.5em},
  main node green/.style={rectangle, draw=mygreen, text width=4.1em,  minimum width = 5.5em, text centered, rounded corners, minimum height=1.5em},
  main node orange/.style={rectangle, draw=orange, text width=4.1em, minimum width = 5.5em, text centered, rounded corners, minimum height=1.5em},
  main node red/.style={rectangle, draw=myred, thick, text width=4.4em,  minimum width = 5.5em, text centered, rounded corners, minimum height=1.5em}
  ]

    \node (generator) {Operations};
\node[main node blue] (4) [below=0.2cm of generator] {};
\node[main node blue] (1) [left=0.2cm of 4] {};
\node[main node green] (2) [below=0.2cm of 1] {$k, p_{jk}, s_{jk}$};
\node[main node red] (A) [left=0.2cm of 5] {$\vphantom{\begin{array}{c}   1\\ 1\\ 1\\ 1\\ 1\\ 1 \end{array}}$};

\node[main node orange] (3) [below=0.2cm of 2] {};
\node[main node green] (5) [right=0.2cm of 2] {$h, p_{jh}, -1$};
\node[main node green] (6) [right=0.2cm of 5] {};


\node (left-paren) [left=0.1cm of 2] {$\left[\vphantom{\begin{array}{c}     1\\ 1\\ 1\\ 1\\ 1\\ 1\\ 1\\ 1\\ 1 \end{array}}\right.$};
\node (right-paren) [right=0.1cm of 6] {$\left.\vphantom{\begin{array}{c}   1\\ 1\\ 1\\ 1\\ 1\\ 1\\ 1\\ 1\\ 1 \end{array}}\right]$};

\node (left-paren1) [left=0cm of 2] {$\left[\vphantom{\begin{array}{c}     1\\ 1 \end{array}}\right.$};
\node (right-paren1) [right=0cm of 6] {$\left.\vphantom{\begin{array}{c}     1\\ 1 \end{array}}\right]$};
\node (left-paren2) [left=0cm of 1] {$\left[\vphantom{\begin{array}{c}     1\\ 1 \end{array}}\right.$};
\node (right-paren2) [right=0cm of 4] {$\left.\vphantom{\begin{array}{c}     1\\ 1 \end{array}}\right]$};
\node (left-paren3) [left=0cm of 3] {$\left[\vphantom{\begin{array}{c}     1\\ 1 \end{array}}\right.$};
\node (right-paren3) [right=0cm of 3] {$\left.\vphantom{\begin{array}{c}     1\\ 1 \end{array}}\right]$};


 \node (1b) [left=0.3cm of 1] {$J_1$};
\node (2b) [left=0.3cm of 2] {$J_j$};
\node (3b) [left=0.3cm of 3] {$J_n$};
 \node (00) [left=0.72cm of 2] {Jobs};
\draw [-, dotted, thick] (1.south) to  (2.north);
\draw [-, dotted, thick] (2.south) to  (3.north);
\draw [-, dotted, thick] (1b.south) to  (2b.north);
\draw [-, dotted, thick] (2b.south) to  ($(3b.north)+(0.18mm,0)$);

\end{tikzpicture}
    \caption{
    Action space $\mathcal{A}_{s_t}$.}
  \label{fig:next_state0}
\end{subfigure}%
\hfill
\\[3ex]
\begin{subfigure}{0.48\textwidth}
  \centering 
  \begin{tikzpicture}[->,>=stealth',auto,node distance=1.42cm,
  semithick,
  main node blue/.style={rectangle, draw=myblue,text width=4.1em, minimum width = 5.5em, text centered, rounded corners, minimum height=1.5em},
  main node green/.style={rectangle, draw=mygreen, text width=4.1em, minimum width = 5.5em, text centered, rounded corners, minimum height=1.5em},
  main node orange/.style={rectangle, draw=orange, text width=4.1em, minimum width = 5.5em, text centered, rounded corners, minimum height=1.5em},
  main node red/.style={rectangle, draw=myred, thick, text width=4.4em, minimum width = 5.5em, text centered, rounded corners, minimum height=1.8em}
  ]

\node (generator) {Operations};
\node[main node blue] (4) [below=0.2cm of generator] {};
\node[main node blue] (1) [left=0.2cm of 4] {};
\node[main node green] (2) [below=0.2cm of 1] {$k, p_{jk}, s_{jk}$};
\node[main node red] (A) [left=0.2cm of 5] {$\vphantom{ciao}$};

\node[main node orange] (3) [below=0.2cm of 2] {};
\node[main node green] (5) [right=0.2cm of 2] {$h, p_{jh}, -1$};
\node[main node green] (6) [right=0.2cm of 5] {};


\node (left-paren) [left=0.1cm of 2] {$\left[\vphantom{\begin{array}{c}     1\\ 1\\ 1\\ 1\\ 1\\ 1\\ 1\\ 1\\ 1 \end{array}}\right.$};
\node (right-paren) [right=0.1cm of 6] {$\left.\vphantom{\begin{array}{c}   1\\ 1\\ 1\\ 1\\ 1\\ 1\\ 1\\ 1\\ 1 \end{array}}\right]$};

\node (left-paren1) [left=0cm of 2] {$\left[\vphantom{\begin{array}{c}     1\\ 1 \end{array}}\right.$};
\node (right-paren1) [right=0cm of 6] {$\left.\vphantom{\begin{array}{c}     1\\ 1 \end{array}}\right]$};
\node (left-paren2) [left=0cm of 1] {$\left[\vphantom{\begin{array}{c}     1\\ 1 \end{array}}\right.$};
\node (right-paren2) [right=0cm of 4] {$\left.\vphantom{\begin{array}{c}     1\\ 1 \end{array}}\right]$};
\node (left-paren3) [left=0cm of 3] {$\left[\vphantom{\begin{array}{c}     1\\ 1 \end{array}}\right.$};
\node (right-paren3) [right=0cm of 3] {$\left.\vphantom{\begin{array}{c}     1\\ 1 \end{array}}\right]$};


 \node (1b) [left=0.3cm of 1] {$J_1$};
\node (2b) [left=0.3cm of 2] {$J_j$};
\node (3b) [left=0.3cm of 3] {$J_n$};
 \node (00) [left=0.72cm of 2] {Jobs};
\draw [-, dotted, thick] (1.south) to  (2.north);
\draw [-, dotted, thick] (2.south) to  (3.north);
\draw [-, dotted, thick] (1b.south) to  (2b.north);
\draw [-, dotted, thick] (2b.south) to  ($(3b.north)+(0.18mm,0)$);

\end{tikzpicture}
    \caption{Action selection $a_t$.}
  \label{fig:action}
\end{subfigure}
\hfill
\begin{subfigure}{0.48\textwidth}\centering 
  \begin{tikzpicture}[->,>=stealth',auto,node distance=1.42cm,
  semithick,
  main node blue/.style={rectangle, draw=myblue,text width=4.1em, minimum width = 5.5em, text centered, rounded corners, minimum height=1.5em},
  main node green/.style={rectangle, draw=mygreen, text width=4.1em, minimum width = 5.5em, text centered, rounded corners, minimum height=1.5em},
  main node orange/.style={rectangle, draw=orange, text width=4.1em, minimum width = 5.5em, text centered, rounded corners, minimum height=1.5em},
  main node faded/.style={rectangle, draw=orange!0, text width=4.1em, minimum width = 5.5em, text centered, rounded corners, minimum height=1.5em}
  ]

\node (generator) {Operations};
\node[main node blue] (4) [below=0.2cm of generator] {};
\node[main node blue] (1) [left=0.2cm of 4] {};
\node[main node green] (2) [below=0.2cm of 1] {$h, p_{jh}, s_{jh}$};
\node[main node orange] (3) [below=0.2cm of 2] {};
\node[main node green] (5) [right=0.2cm of 2] {};
\node[main node faded] (6) [right=0.2cm of 5] {};


\node (left-paren) [left=0.1cm of 2] {$\left[\vphantom{\begin{array}{c}     1\\ 1\\ 1\\ 1\\ 1\\ 1\\ 1\\ 1\\ 1 \end{array}}\right.$};
\node (right-paren) [right=0.1cm of 6] {$\left.\vphantom{\begin{array}{c}   1\\ 1\\ 1\\ 1\\ 1\\ 1\\ 1\\ 1\\ 1 \end{array}}\right]$};

\node (left-paren1) [left=0cm of 2] {$\left[\vphantom{\begin{array}{c}     1\\ 1 \end{array}}\right.$};
\node (right-paren1) [right=0cm of 5] {$\left.\vphantom{\begin{array}{c}     1\\ 1 \end{array}}\right]$};
\node (left-paren2) [left=0cm of 1] {$\left[\vphantom{\begin{array}{c}     1\\ 1 \end{array}}\right.$};
\node (right-paren2) [right=0cm of 4] {$\left.\vphantom{\begin{array}{c}     1\\ 1 \end{array}}\right]$};
\node (left-paren3) [left=0cm of 3] {$\left[\vphantom{\begin{array}{c}     1\\ 1 \end{array}}\right.$};
\node (right-paren3) [right=0cm of 3] {$\left.\vphantom{\begin{array}{c}     1\\ 1 \end{array}}\right]$};


 \node (1b) [left=0.3cm of 1] {$J_1$};
\node (2b) [left=0.3cm of 2] {$J_j$};
\node (3b) [left=0.3cm of 3] {$J_n$};
 \node (00) [left=0.72cm of 2] {Jobs};
\draw [-, dotted, thick] (1.south) to  (2.north);
\draw [-, dotted, thick] (2.south) to  (3.north);
\draw [-, dotted, thick] (1b.south) to  (2b.north);
\draw [-, dotted, thick] (2b.south) to  ($(3b.north)+(0.18mm,0)$);

\end{tikzpicture}
    \caption{Next state representation $s_{t+1}$.}
  \label{fig:next_state}
\end{subfigure}
\caption{Graphical representation of the environment and its components.}
\label{fig:iteration}
\end{figure}

It turns out  this process can be formulated as a Markov Decision process (MDP). With a notation similar  to the one used to the describe RL in section \ref{sec:preliminaries_and_notation}, we characterize a finite MDP using the tuple  $(\mathcal{S},\mathcal{A},\mathcal{R},P)$, where $\mathcal{S}$ is the set of states, $\mathcal{A}$  the set of actions, ${R}: \mathcal{S}\times\mathcal{A}\rightarrow\mathbb{R}$ a reward function and $P:\mathcal{S}\times\mathcal{A}\rightarrow\mathcal{S}$ the transition function. The \textbf{state} ${s_t}\in \cal S$ captures all the relevant information regarding the current iteration, in order to respect the {Markov Property} and to have a  fully observable system state.
The data structure adopted is a list of jobs, with each job being a list of tasks (or operations) $(j,k)$, with $j \in \cal J$ being the job and $k \in \cal M$ the machine, as illustrated in Figure \ref{fig:state}. A processing time $p_{jk}$ is associated with each operation $(j,k)$, representing the time spent by the job $j$ to complete the task on the machine $k$. If a task $(j,k)$ is the first available operation for a  job, then we take into account  the \textit{earliest possible starting time} $s_{jk}$, defined as 
$s_{jk} = \max(C_{jh},C_{ik})$, where  $(j,h)~\prec~(j,k)$, 
$(i,k)$ is  the last operation scheduled using the machine $k$,
and  $C_{jh}$ and $C_{ik}$ are the corresponding completion times.
At each decision step $t$, an \textbf{action} $a_t$, representing the allocation of an operation at a certain starting time, is taken and the corresponding task is removed from the job and put into a list of scheduled operations.  
The set ${\mathcal{A}_{s_t}}$ contains all the available actions at state $s_t$, i.e. the operations still to be scheduled, with no prior operations or whose prior operations have been already scheduled, as illustrated in Figures   \ref{fig:next_state0} and \ref{fig:action}. The \textbf{next state} ${s_{t+1}}$ is the state $s_{t}$ without the operation scheduled at decision step $t$ (see  Figure \ref{fig:next_state}). Due to the deterministic  nature of the problem, the state transition is deterministic and given a state-action pair $(s_t,a_t)$, the next state $s_{t+1}$ is uniquely determined. In order to be coherent with the RL classical notation, the \textbf{reward} ${r_{t+1}}$  will be the negative contribution of the selected operation to the current makespan $C_{\max_t}$.
If $(j,k)$ is the operation scheduled at decision step $t$, then
 
 $$
 r_{t+1} = 
    \begin{cases}
       -(C_{jk}-C_{\max_t}),& \text{if } C_{jk} > C_{\max_t}\\
       0,& \text{otherwise} 
    \end{cases}
 $$
 
 In this way, it holds  $C_{\max_k} = \displaystyle -\sum_{t=1}^{k} r_t$, and $C_{\max} = \displaystyle -\sum_{t=1}^{T} r_t$.

\section{Agents deep models}\label{sec:agents_deep_models}

In this section, we describe the neural network infrastructures behind our algorithm. In particular, since we are adopting an actor-critic method, we need to specify two networks: the first, referred to as the actor, estimates the policy $\pi_\theta(\cdot,s_t)$, while the second,  the critic, provides an estimation of the state-value function $V_\omega(s_t)$. Our aim is to obtain a learned algorithm that is flexible concerning the size of the JSSP instance, both in the number of jobs and machines.
According to the state defined in Section \ref{sec:jssp_as_mdp}, { the state is represented by a list of tasks for each job, translatind to a  }data structure characterized by a list, variable in size,  of lists, with different lengths the one to the other.
For these reasons, the actor and the critic networks make use of long short-term memory networks (LSTMs), which are suitable for processing variable-length sequences. {LSTMs process the information in sequence-structured input, by taking its elements one at a the time, eventually propagating the information using special arcs.} LSTMs belong to the class of recurrent neural networks, and they deal with the problem of the vanishing gradient thanks to these special arcs called self-loops in the hidden layer, as described in \cite{graves2012long}.
The self-loops are controlled by the network, which is able to adjust the information flow.
They take sequences as input, returning same-length sequences of embeddings.
\medskip

\begin{figure}
\centering 
\begin{subfigure}{0.4\textwidth}
\centering 
  \begin{tikzpicture}[->,>=stealth',auto,node distance=1.42cm,
  thick,
  main node blue/.style={rectangle, draw=myblue,text width=2.6em, text centered, rounded corners, minimum height=6.1em, minimum width = 3em},
  main node green/.style={rectangle, draw=mygreen, text width=4.1em, text centered, rounded corners, minimum height=2.8em, minimum width = 13.5em},
  main node orange/.style={rectangle, draw=orange, text centered, text width=8.1em, minimum height=8em, minimum width = 13.5em}
  ]
    \footnotesize
    
\node (generator) {State};
\draw[-] [decorate,semithick,decoration={brace,amplitude=5pt,mirror,raise=4ex}]
  (-2,0.2) -- (2,0.2) node[midway,yshift=+3em]{};
\node[main node blue] (2) [below=0.5cm of generator] {LSTM 1\\\ \\$J_j$};
\node[main node blue] (1) [left of=2] {LSTM 1\\\ \\$J_1$};
\node[main node blue] (3) [right of=2] {LSTM 1\\\ \\$J_n$};

\node[main node green] (4) [below=0.2cm of 2] {LSTM 2};
\node[main node orange] (5) [below=0.8cm of 4] {Action Masking \\+ \\Softmax $\sigma$};

\node (out)  [below = 0.8cm of 5] {$\hat{\pi}(\cdot,s)$};




 
\draw [-, dotted, thick] ($(1.east) + (1mm,0)$) to  ($(2.west) - (1mm,0)$);
\draw [-, dotted, thick] ($(2.east)+ (1mm,0)$) to  ($(3.west)- (1mm,0)$);
\draw [->, semithick] ($(4.south)-(0,1.5mm)$) node[ below] {$\quad\qquad\qquad  y\in\mathbb{R}^{|\mathcal{J}|}$}   to  ($(5.north)+(0,1.5mm)$);
\draw [->, semithick] ($(5.south)-(0,1.5mm)$) to  ($(out.north)+(0,1mm)$);


\end{tikzpicture}
    \caption{Actor model.} \label{fig:Actor_model}
\end{subfigure}%
\hfill
\begin{subfigure}{0.4\textwidth}
  \centering 
  \begin{tikzpicture}[->,>=stealth',auto,node distance=1.42cm,
  thick,
  main node blue/.style={rectangle, draw=myblue, text width=2.6em, text centered, rounded corners, minimum height=6.1em, minimum width = 3em},
  main node green/.style={rectangle, draw=mygreen, text width=4.1em, text centered, rounded corners, minimum height=2.8em, minimum width = 13.5em},
  main node orange/.style={signal, draw=myred, signal to=south, text centered, align=center, text width=8.1em, minimum height=0em, minimum width = 13.5em}
  ]
   \footnotesize
   
\node (generator) {State};
\draw[-] [decorate, semithick, decoration={brace,amplitude=5pt,mirror,raise=4ex}]
  (-2,0.2) -- (2,0.2) node[midway,yshift=+3em]{};
\node[main node blue] (2) [below=0.5cm of generator] {LSTM 1\\\ \\$J_j$};
\node[main node blue] (1) [left of=2] {LSTM 1\\\ \\$J_1$};
\node[main node blue] (3) [right of=2] {LSTM 1\\\ \\$J_n$};

\node[main node green] (4) [below=0.2cm of 2] {LSTM 2};
\node [rectangle, draw=green!0, text width=4.1em, text centered, rounded corners, minimum height=2.8em, minimum width = 4.1em] (6) [below=1.2cm of 4, xshift=-3.5] {Feedforward\\ \ \ \ (FFN)};
\node[main node orange] (5) [below=0.8cm of 4] {};

\node (out)  [below = 0.8cm of 5] {$\hat{V}(s)$};




 
\draw [-, dotted, thick] ($(1.east) + (1mm,0)$) to  ($(2.west) - (1mm,0)$);
\draw [-, dotted, thick] ($(2.east)+ (1mm,0)$) to  ($(3.west)- (1mm,0)$);
\draw [->, semithick] ($(4.south)-(0,1.5mm)$) node[ below] {$\quad\qquad\qquad z\in\mathbb{R}^{h}$}   to  ($(5.north)+(0,1.5mm)$);
\draw [->, semithick] ($(5.south)-(0,1.5mm)$) to  ($(out.north)+(0,1mm)$);


\end{tikzpicture}
    \caption{Critic model.}\label{fig:Critic_model}
\end{subfigure}

\caption{Deep neural infrastructures used.}
\label{fig:iteration1}
\end{figure}
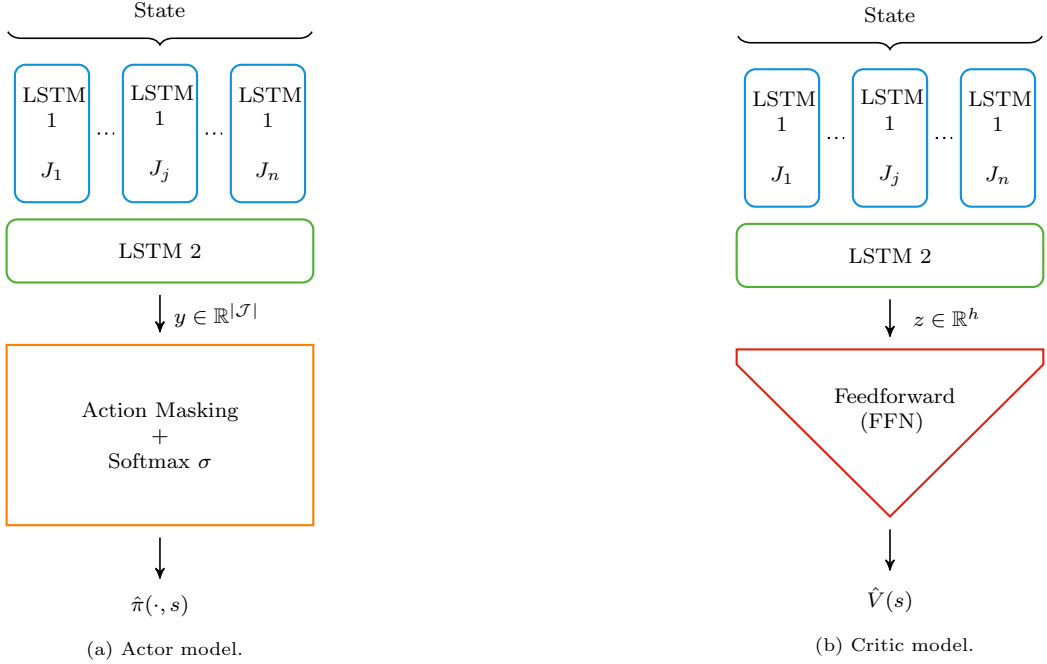

The actor model is composed of two concatenated LSTMs. The first takes as input the state,  producing an embedding for each operation. We then consider only the embedding related to the last operation of each job, since operations in a job are chronologically connected and the last element can be seen as a compressed representation of the whole job. After that, a list with size  $|\mathcal{J}|$  is obtained from the embeddings and passed the second LSTM. This network combines the jobs information as a sequence of $|\mathcal{J}|$ embeddings. Each component is collapsed to a scalar, obtaining a vector $y\in\mathbb{R}^{|\mathcal{J}|}$. {This trick allows us to compress the information and dealwith the two levels of variability of the data structure.}
Finally, we use action masking and a softmax function, applying the mask vector $M$, whose components are Boolean values filtering out the invalid actions at each state. The softmax is used to transform the embedding $y$ into a probability distribution through the formula 
$$ \sigma(y,M)_i = \frac{e^{y_i}M_i}{\sum_{j=1}^{|\mathcal{J}|}e^{y_j}M_j}, \quad i=1,\dots,|\mathcal{J}|$$
The actor network is represented in Figure \ref{fig:Actor_model}.

The critic model is composed of a double LSTM as in the actor, followed by a deep feed-forward neural network (FFN). 
After being fed with the output by the first LSTM, the second LSTM returns a sequence of embeddings with cardinality $|\mathcal{J}|$. Then, the vectors are summed up, obtaining the vector $z\in \setR^h$, where $h$ is the hidden size of the second LSTM. Finally, the FFN processes $z$, converging to a scalar. The FFN has three fully connected hidden layers with a decreasing number of neurons, applying the ReLU as an activation function until the last layer, which is linear.
The critic network is illustrated in Figure \ref{fig:Critic_model}.

Before feeding the two neural models, the input sequences must be padded to level equally their lengths within a single mini-batch.

\section{Computational Experience}\label{sec:computational_experience}

The computational experience exploits both the training and the testing phase. In the former, we give insights on the computational burden required for training, whereas in the latter, we evaluate our algorithm against:
 \begin{itemize}
     \item a branch-and-cut algorithm implemented in the commercial optimization solver CPLEX v12.10, \cite{cplex}, Section \ref{sec:cplex},
      \item an adaptive scheduling algorithm based on the mixed disjunctive graph model, proposed in \cite{Sotskov2013Adaptive} and \cite{Sotskov1996Adaptive}, Section \ref{sec:adaptive},
     \item 17 well-known scheduling heuristics, where at each step an operation is scheduled according to a priority dispatching rule (PDR), \cite{Panwalkar1977ASurvey}, Section \ref{sec:priority_rules},
   \item      { the Deep Reinforcement learning approach based on Graph Neural Networks proposed in \cite{zhang2020learning}, Section \ref{sec:zhang}.}
 \end{itemize}
Our algorithm uses PyTorch v1.8, \cite{pytorch}, working on a Windows server with a11th Gen Intel(R) Core(TM) i7-11800H CPU and a single NVIDIA GeForce RTX 3060 GPU.
In the following, we will use the notation $Jobs \times Machines$ to specify the  dimension of the JSSPs used.

   {The source code and the data of the experiments are available at:

\href{https://github.com/GiorgioGrani/JSSP_actor-critic_Agasucci_Monaci_Grani}{github.com/GiorgioGrani/JSSP\_actor-critic\_Agasucci\_Monaci\_Grani}}

\subsection{Training process}

In the training process, we used a dataset composed   {of} randomly generated JSSP instances, with various, yet small, sizes. We used five classes of JSSP problems,
from  $(8\times6)$ up to  $(20\times15)$. The processing times are drawn from a Gaussian distribution $\mathcal{N}(\mu=100,\sigma=10)$. There is a total of $1000$ instances, having $200$ elements per class. Table \ref{tab:training_sets} summarizes the information regarding the training set.

\begin{table}
	\centering
	\footnotesize
			 \renewcommand\arraystretch{1.5} 
			\begin{tabular}{cccc}
\hline
Jobs $\times$ Machines & $8\times 6$,\ $10\times9,\ 15\times10,\ 17\times13,\ 20\times15$\\

\# Instances & $200,\ 200,\ 200,\ 200,\ 200$\\

{Gaussian} processing times dist. & $p_{jk}\sim\mathcal{N}(\mu = 100,\sigma = 10)$\\
\hline
\end{tabular}
    \smallskip
	\caption{Training set.}
	\label{tab:training_sets}
\end{table}

During the learning process, the agent has to balance between exploration, i.e. trying new actions, and exploitation of the past experience. 
To prevent the agent from getting stuck in bad regions, we  inject the state-space exploration  introducing the  probability $\epsilon$ to  perform a random action.
The value of $\epsilon$ is updated according to a step function starting from $\epsilon=0.20$ and gradually decreases until, at around the 70\% of the total number of episodes, goes down to $\epsilon=0$. We ran the training over $5000$ episodes, with $N=10$ roll-outs per episode.

The actor-critic network configurations and training settings are reported in the Appendix.
The plot in Figure \ref{fig:value_loss} reports the state-value function loss, showing both the single value and the moving average over $100$ episodes. The image indicates a tendency in the reduction of the loss and the variability.

\begin{figure}
    \centering
    \tiny 
    \resizebox{0.65\textwidth}{!}{
    \input{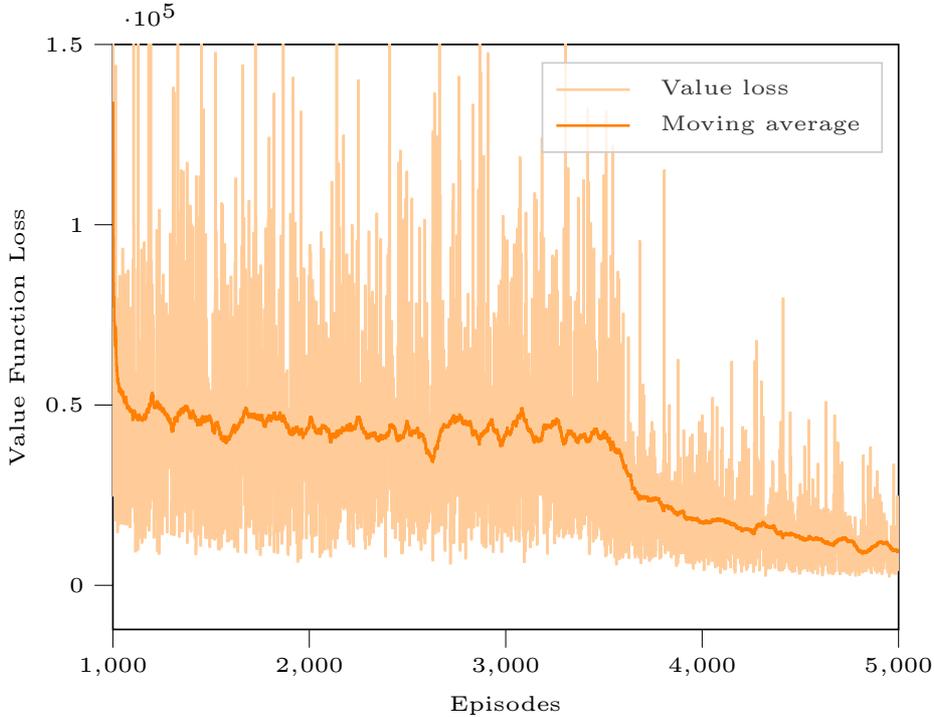}
    }
    \caption{State-value function loss during training process}
    \label{fig:value_loss}
\end{figure}

Since the total reward measure varies from instance to instance, we define a relative objective value gap, $\phi_k$, which is the ratio between the total reward obtained in the $k$-th learning iteration for the instance $i(k)$, and the best solution found trough the all training for the same instance $i(k)$.

$$ \phi_k= \frac{R_k}{\min \left\{ R_h \ : \ h = 1,\dots, 5000\ \wedge\ i(h) = i(k) \right\}},  $$

\begin{figure}
\centering 
\footnotesize
\begin{subfigure}{0.5\textwidth}
\centering 
  \input{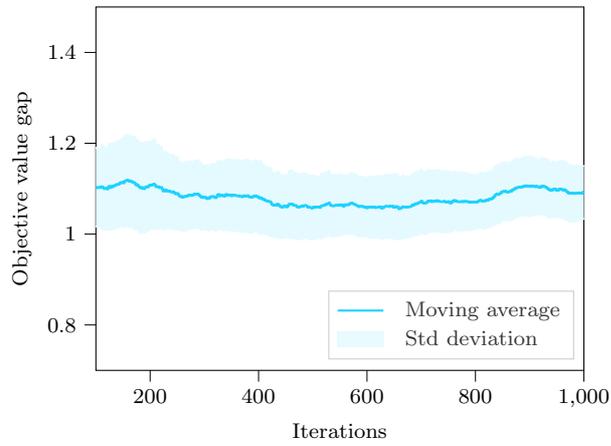}
    \caption{Size $8\times 6$.}
    \label{fig:reward8x6}
\end{subfigure}%
\hfill
\begin{subfigure}{0.5\textwidth}
  \centering 
  \input{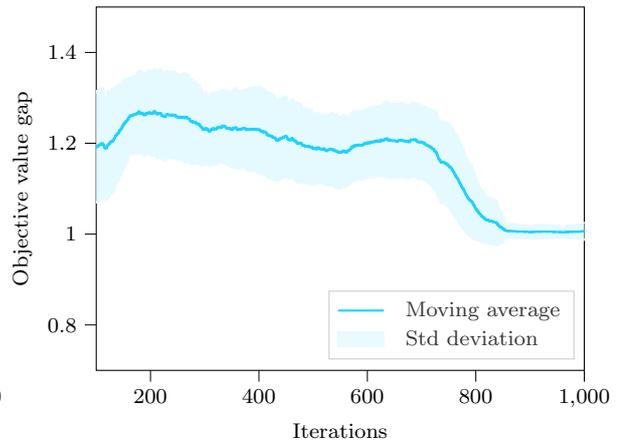}
    \caption{Size $10\times 9$.}
    \label{fig:reward10x9}
\end{subfigure}
\medskip

\begin{subfigure}{0.5\textwidth}
\centering 
  \input{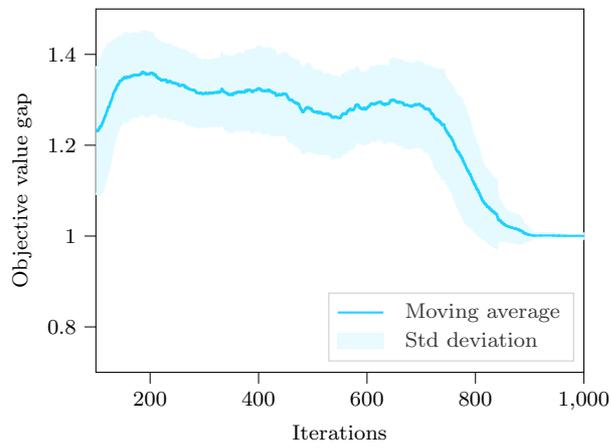}
    \caption{Size $15\times 10$.}
    \label{fig:reward15x10}
\end{subfigure}%
\hfill
\begin{subfigure}{0.5\textwidth}
  \centering 
  \input{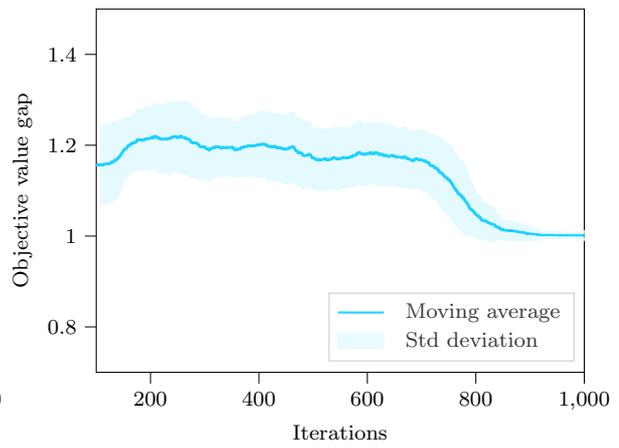}
    \caption{Size $17\times 13$.}
    \label{fig:reward17x13}
\end{subfigure}

\medskip

\begin{subfigure}{0.5\textwidth}
\centering 
  \input{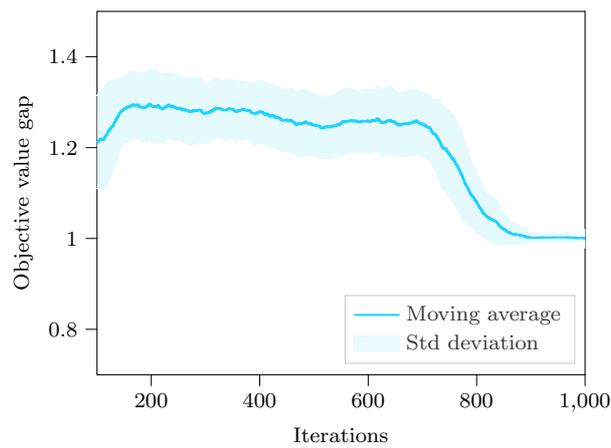}
    \caption{Size $20\times 15$.}
    \label{fig:reward20x15}
\end{subfigure}%

\caption{Objective value gap during the training process.}
\label{fig:rewards.}
\end{figure}

The graphs in Figure \ref{fig:rewards.} illustrate the moving average and the confidence interval for $\phi_k$, divided by classes.
There is a decreasing trend with a jump around the 70\% of the iterations, in accordance with the random choice probability $\epsilon$.

   {Regarding the computational time required to perform the training, it takes around $8$ hours on a Windows machine with a 11th Gen Intel(R) Core(TM) i7-11800H CPU and a single NVIDIA GeForce RTX 3060 GPU.}

\subsection{Test process}

The test phase has been conducted to demonstrate the ability of our algorithm to generalize over larger instances, generalize over different distributions, and maintain computational efficiency. To this aim, we tested the performances of our approach on two sets with five JSSP classes each, going from $(30\times25)$ to $(50\times45)$, with 100 instances per class. In the first dataset, the processing times are derived from a Gaussian distribution $\mathcal{N}(\mu=100,\sigma=10)$, while in the second they are drawn from a Poisson $\mathcal{P}(\lambda=100)$. Table \ref{tab:test_sets} summarizes the specifics of the test sets.

\begin{table}
	\centering
	\footnotesize
			 \renewcommand\arraystretch{1.5} 
			\begin{tabular}{cccc}
\hline
Jobs $\times$ Machines & $30 \times 25$,\ $35\times30,\ 40\times35,\ 45\times40,\ 50\times45$\\
 
\# Instances & $100,\ 100,\ 100,\ 100,\ 100$\\
 
{Gaussian} processing times dist. & $p_{jk}\sim\mathcal{N}(\mu = 100,\sigma = 10)$\\

{Poisson} processing times dist.& $p_{jk}\sim\mathcal{P}(\lambda = 100)$\\
\hline
\end{tabular}
    \smallskip
	\caption{Test sets: {Gaussian} and {Poisson sets}.}
	\label{tab:test_sets}
\end{table}

    {In Sections \ref{sec:cplex}, \ref{sec:priority_rules}, and \ref{sec:zhang}, we introduced some comparisons on $50$ Taillard's benchmark instances taken from \cite{taillard1993benchmarks}. The size of these instances are resumed in Table \ref{tab:taillard}, and we will refer to them as the Taillard's Benchmark set.}

    {In Section \ref{sec:zhang}, we will introduce an additional set of generated Taillard instances, and we redirect the description of these sets to the dedicated Section.}

\begin{table}
	\centering
	\footnotesize
			 \renewcommand\arraystretch{1.5} 
			   {\begin{tabular}{ccc}
\hline
Jobs $\times$ Machines & Instance IDs & $\#$instances\\ \hline
$15 \times 15$ & Ta01-10 & 10 \\
$20 \times 15$ & Ta01-10 & 10 \\
$20 \times 20$ & Ta01-10 & 10 \\
$30 \times 15$ & Ta01-10 & 10 \\
$30 \times 20$ & Ta01-10 & 10 \\ \hline
\end{tabular}}
    \smallskip
	\caption{Taillard's benchmark set.}
	\label{tab:taillard}
\end{table}

\subsubsection{Comparison with the branch-and-cut algorithm implemented in CPLEX}\label{sec:cplex} 

We tested out our algorithm against the branch-and-cut algorithm
implemented in the known commercial solver CPLEX v12.10, \cite{cplex}, conducting two types of analyses:
\begin{itemize}
    \item {\textit{Objective value analysis}}, where we study the objective value, i.e. the makespan $C_{max}$, reached by the algorithm in CPLEX when its time limit is locked to be no more than our timing.
    
    \item \textit{Computational time analysis}, where we compare our timing with the one required for the algorithm in CPLEX to reach a solution as good as ours in terms of objective value, i.e. the makespan $C_{max}$.

\end{itemize}

The results of the objective value analysis on the Gaussian and Poisson sets are presented in Tables \ref{tab:res_gauss_obj} and \ref{tab:res_poiss_obj}, respectively. These show statistics on objective value, pruned from outliers.
We use the statistic $\rho$ to express the percentage of average  improvement of our algorithm {in terms of the objective value $C_{\max}$. Let $\rho_i$ be the improvement on the single $i$-th instance}
$$
\rho_i = \frac{C_{\max_{CP_i}}-C_{\max_{RL_i}}}{C_{\max_{CP_i}}},
$$
where $C_{\max_X}$ denotes the makespan of the algorithm $X$ and
X can equal $RL$ (the proposed RL method) or $CP$ (short for the branch-and-cut procedure implemented in BC-CPLEX).

Our approach outperforms the branch and cut algorithm implemented in CPLEX in terms of objective value. For the $(40\times35)$  class in the Gaussian set, our approach returns a makespan on average 91\% better than the branch-and-cut one, and 87\% in the {Poisson set} for the same class. It is worth mentioning that, for no instance, the algorithm implemented in CPLEX terminated before the time limit. The values of $\rho$
tend to improve for larger classes, since instances become more and more complex for the deterministic solver.

\begin{table}[H]
\centering
\footnotesize
			 \scalebox{0.92}{
			 \renewcommand\arraystretch{1.5} 
\begin{tabular}{ccccccccccc}
\cline{1-11}
\multicolumn{1}{c}{J $\times$ M} &
  \multicolumn{1}{c}{$\#$instances} &
  \multicolumn{2}{c}{Mean} &
  \multicolumn{2}{c}{Std dev} &
  \multicolumn{2}{c}{Max} &
  \multicolumn{2}{c}{Min} &
  Avg $\rho$ (\%) \\ \hline
 &
   &
  \multicolumn{1}{c}{RL} &
  \multicolumn{1}{c}{BC-CPLEX} &
  \multicolumn{1}{c}{RL} &
  \multicolumn{1}{c}{BC-CPLEX} &
  \multicolumn{1}{c}{RL} &
  \multicolumn{1}{c}{BC-CPLEX} &
  \multicolumn{1}{c}{RL} &
  \multicolumn{1}{c}{BC-CPLEX} &
  \multicolumn{1}{c}{} \\ \cline{3-10}
30 $\times$ 25 & 100 & \textbf{4785.3} & 29505.2 & 111.7 & 15261.5 & 5107.3 & 40408.4  & 4579.7 & 4307.1  & \textbf{63.1} \\
35 $\times$ 30 & 100 & \textbf{5462.6} & 53883.9 & 115.1 & 11398.5 & 5702.0 & 60350.1  & 5180.4 & 5309.8  & \textbf{85.8} \\
40 $\times$ 35 & 100 & \textbf{5989.5} & 68449.6 & 134.8 & 4592.1  & 6349.8 & 74838.7  & 5715.4 & 62320.4 & \textbf{91.2} \\
45 $\times$ 40 & 100 & \textbf{6931.4} & 66026.5 & 130.2 & 38853.8 & 7304.9 & 99942.2  & 6599.8 & 6376.8  & \textbf{63.9} \\
50 $\times$ 45 & 100 & \textbf{8663.5} & 93547.3 & 147.4 & 44467.5 & 9015.6 & 125757.1 & 8198.9 & 7827.4  & \textbf{71.6} \\ 
\end{tabular}}
    \smallskip
	\caption{{RL vs the branch-and-cut implemented in CPLEX (BC-CPLEX) on the {Gaussian set}}: comparison on the makespan $C_{max}$.}
	\label{tab:res_gauss_obj}
\end{table}

\begin{table}[H]
\centering
\footnotesize
\scalebox{0.92}{
\renewcommand\arraystretch{1.5} 
\begin{tabular}{ccccccccccc}
\cline{1-11}
\multicolumn{1}{c}{J $\times$ M} &
  \multicolumn{1}{c}{$\#$instances} &
  \multicolumn{2}{c}{Mean} &
  \multicolumn{2}{c}{Std dev} &
  \multicolumn{2}{c}{Max} &
  \multicolumn{2}{c}{Min} &
  Avg $\rho$ $(\%)$ \\ \hline
 &
   &
  \multicolumn{1}{c}{RL} &
  \multicolumn{1}{c}{BC-CPLEX} &
  \multicolumn{1}{c}{RL} &
  \multicolumn{1}{c}{BC-CPLEX} &
  \multicolumn{1}{c}{RL} &
  \multicolumn{1}{c}{BC-CPLEX} &
  \multicolumn{1}{c}{RL} &
  \multicolumn{1}{c}{BC-CPLEX} &
  \multicolumn{1}{c}{} \\ \cline{3-10}
30 $\times$ 25 & 100 & \textbf{4167.1} & 28079.9  & 93.7  & 12378.8 & 4404.0 & 36509  & 3939.0 & 4058.0 & \textbf{69.8} \\
35 $\times$ 30 & 100 & \textbf{5249.4} & 43213.3  & 136.0 & 20358.0 & 5566.0 & 55073.0  & 4963.0 & 4609.0 & \textbf{69.8} \\
40 $\times$ 35 & 100 & \textbf{5947.5} & 72325.4  & 111.9 & 16863.0 & 6194.0 & 78182.0  & 5735.0 & 5777.0 & \textbf{86.8} \\
45 $\times$ 40 & 100 & \textbf{7815.5} & 71007.8 & 182.6 & 34066.1 & 8263.0 & 96754.0  & 7453.0 & 6702.0 & \textbf{68.9} \\
50 $\times$ 45 & 100 & \textbf{8032.1} & 103049.2 & 181.9 & 40447.6 & 8410.0 & 129359.0 & 7697.0 & 7571.0 & \textbf{78.6} \\ 
\end{tabular}}
    \smallskip
	\caption{{RL vs branch-and-cut implemented in CPLEX (BC-CPLEX) on the {Poisson set}}: comparison on the makespan $C_{max}$.}
	\label{tab:res_poiss_obj}
\end{table}

The results of the computational time analysis on the Gaussian and Poisson sets are summarized in Tables \ref{tab:res_gauss_time} and \ref{tab:res_poiss_time}, respectively. These tables report the statistics on computational time for both algorithms. We use the statistic $\tau$ to express the percentage of average improvement of our algorithm {in terms of the computational time. Let $\tau_i$ be the improvement on the single $i$-th instance}
$$
\tau_i = \frac{time_{CP_i}-time_{RL_i}}{time_{CP_i}},
$$
where ${time}_X$ denotes the computational time of the algorithm $X$ and
$X$ can equal $RL$ (the proposed RL method) or $CP$ (short for the branch-and-cut algorithm
implemented in the known commercial solver BC-CPLEX).

For each class of both sets, our method is faster on average than the optimization solver implemented in CPLEX, up to reach a  78\% improvement for the $(40\times35)$ instances in the Gaussian set, and 67\% for the Poisson set. 
For some classes, the number of instances is less than 100 since the branch-and-cut algorithm
implemented in the known commercial solver CPLEX  exceeded the time limit of two minutes.

\begin{table}[H]
\centering
\footnotesize
			 \renewcommand\arraystretch{1.5} 
\begin{tabular}{ccccccccccc}
\hline
J $\times$ M     & \#instances & \multicolumn{2}{c}{Mean} & \multicolumn{2}{c}{Std dev} & \multicolumn{2}{c}{Max} & \multicolumn{2}{c}{Min} & Avg $\tau$ $(\%)$ \\ \hline
               &             & RL              & BC-CPLEX  & RL          & BC-CPLEX         & RL         & BC-CPLEX      & RL         & BC-CPLEX      &                   \\ \cline{3-10}
30 $\times$ 25 & 100         & \textbf{2.5}    & 4.7    & 0.3         & 2.4           & 3.5        & 14.4       & 1.9        & 1.9        & \textbf{35.8}     \\
35 $\times$ 30 & 100         & \textbf{4.6}    & 13.6   & 0.5         & 6.8           & 5.6        & 36.4       & 3.5        & 3.6        & \textbf{58.2}     \\
40 $\times$ 35 & 99          & \textbf{7.0}    & 39.3   & 1.0         & 16.3          & 10.3       & 79.6       & 5.8        & 12.8       & \textbf{78.3}     \\
45 $\times$ 40 & 100         & \textbf{11.8}   & 34.0   & 1.2         & 22.3          & 14.0       & 99.3       & 9.8        & 5.2        & \textbf{36.4}     \\
50 $\times$ 45 & 92          & \textbf{18.7}   & 51.6   & 1.9         & 34.1          & 24.4       & 117.9      & 15.9       & 9.1        & \textbf{37.1}     \\ 
\end{tabular}
    \smallskip
	\caption{{RL vs the branch-and-cut algorithm
implemented in the known commercial solver CPLEX (BC-CPLEX)} on the {Gaussian set}: comparison on the computational time (seconds).}
	\label{tab:res_gauss_time}
\end{table}

\begin{table}[H]
\centering
\footnotesize
			 \renewcommand\arraystretch{1.5} 
\begin{tabular}{ccccccccccc}
\hline
J $\times$ M     & \#instances & \multicolumn{2}{c}{Mean} & \multicolumn{2}{c}{Std dev} & \multicolumn{2}{c}{Max} & \multicolumn{2}{c}{Min} & Avg $\tau$ $(\%)$    \\ \hline
               &             & RL              & BC-CPLEX  & RL          & BC-CPLEX         & RL         & BC-CPLEX      & RL         & BC-CPLEX      &               \\ \cline{3-10}
30 $\times$ 25 & 100         & \textbf{2.4}    & 5.6    & 0.4         & 2.7           & 4.7        & 13.2       & 1.8        & 1.9        & \textbf{48.3} \\
35 $\times$ 30 & 100         & \textbf{4.4}    & 8.7    & 0.5         & 4.3           & 6.1        & 27.2       & 3.4        & 3.3        & \textbf{39.6} \\
40 $\times$ 35 & 100         & \textbf{7.2}    & 28.4   & 0.9         & 15.1          & 10.7       & 91.8       & 6.0        & 6.9        & \textbf{66.9} \\
45 $\times$ 40 & 100         & \textbf{11.3}   & 30.8   & 1.1         & 20.0          & 14.1       & 79.4       & 9.7        & 5.4        & \textbf{37.5} \\
50 $\times$ 45 & 98          & \textbf{17.1}   & 57.2   & 1.2         & 34.1          & 20.2       & 106.0      & 15.0       & 11.8       & \textbf{54.3} \\ 
\end{tabular}
    \smallskip
	\caption{{RL vs the branch-and-cut algorithm
implemented in the known commercial solver CPLEX} on the {Poisson set}: comparison on the computational time (seconds).}
	\label{tab:res_poiss_time}
\end{table}

\medskip

\begin{figure}
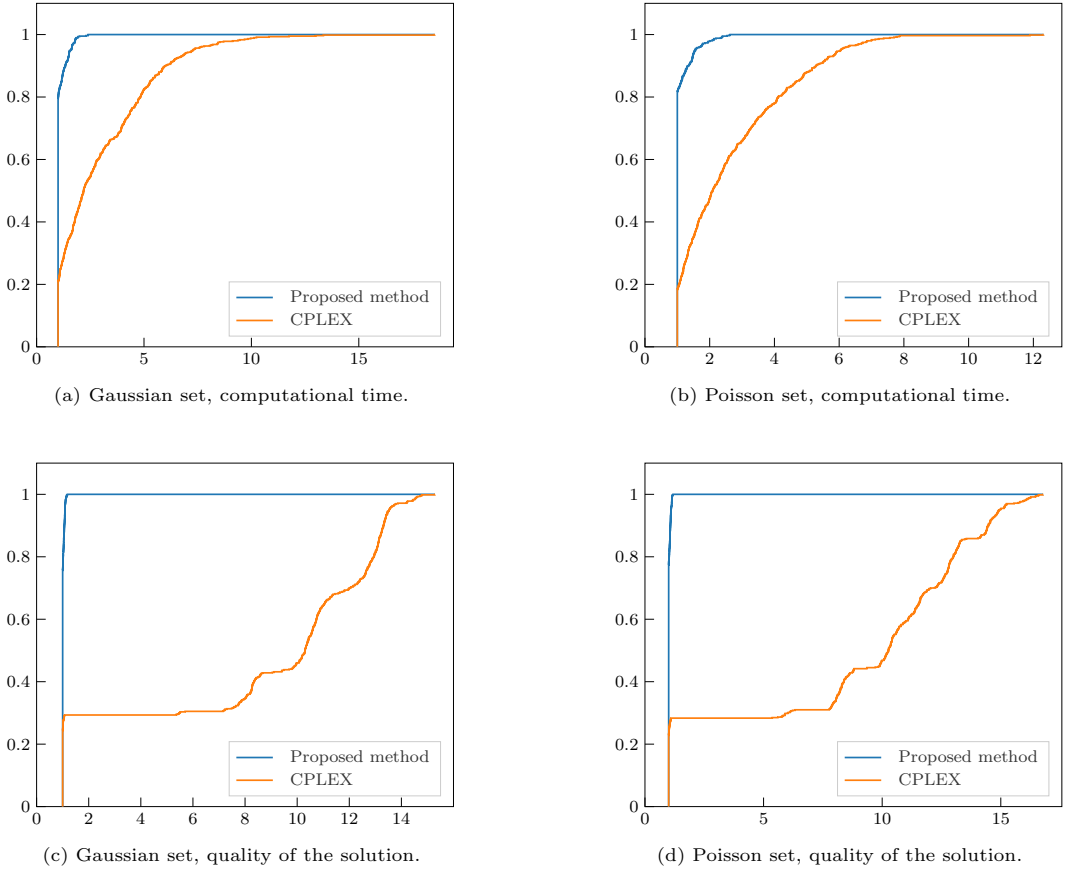

\centering 
\footnotesize
\begin{subfigure}{0.5\textwidth}
\centering 
  \input{gaussian_time}
    \caption{Gaussian set, computational time.}
    \label{fig:gaussiantime}
\end{subfigure}%
\hfill
\begin{subfigure}{0.5\textwidth}
  \centering 
  \input{poisson_time}
    \caption{Poisson set, computational time.}
    \label{fig:poissontime}
\end{subfigure}

\vspace{.7cm}

\begin{subfigure}{0.5\textwidth}
\centering 
  \input{gaussian_obj_value}
    \caption{Gaussian set, objective value.}
    \label{fig:gaussianquality}
\end{subfigure}%
\hfill
\begin{subfigure}{0.5\textwidth}
  \centering 
  \input{poisson_obj_value}
    \caption{Poisson  set, objective value.}
    \label{fig:poissonquality}
\end{subfigure}

\caption{Performance profiles.}
\label{fig:performance_profiles}
\end{figure}

In Figure \ref{fig:performance_profiles}, we report the performance profiles for the Gaussian and the Poisson set, considering the computational time analysis first, and the objective value afterwords.
Performance profiles were initially introduced in \cite{DM2002} as an additional way to compare different methods.
Given a set of algorithms  $\mathcal{A}$ and a set of problems $\mathcal{P}$,
the performance of an algorithm $a \in \cal A$ on a problem $p \in \cal P$ is computed against the best performance obtained by any other method in $\cal A$ on $p$. We consider the ratio
$
\eta_{p,a} = \textit{performance}_{p,a}/\min\{\textit{performance}_{p,a^\prime}: a^\prime \in\mathcal{A}\},
$
where $\textit{performance}_{p,a}$ is the performance obtained on the $p-$th problem by the $a-$th algorithm. In our case, the performance is the computational time at first, and then the objective value. We now consider a cumulative function computing the number of times algorithm $a \in \cal A$ was successful against the others, specifically $\gamma_a(\tau) = |\{p\in \mathcal{P}:\; \eta_{p,a}\leq \tau \}| /|\mathcal{P}|$.
The performance profile is the plot of the function $\gamma_a(\tau)$ for all $a \in \cal A$, varying with $\tau$.

The plots in Figure \ref{fig:performance_profiles} are self-explanatory, since our algorithm outperforms the branch-and-cut algorithm
implemented in the known commercial solver CPLEX in every case. In particular, Figure \ref{fig:gaussiantime} and Figure \ref{fig:poissontime} show that our approach is faster than 80\% of all the times produced by the algorithm implemented in CPLEX, independently from the distribution. For Figure \ref{fig:gaussianquality} and \ref{fig:poissonquality}, the situation is even more accentuated, with our approach beating all the instances almost immediately, whereas the branch-and-cut procedure in CPLEX requires a relatively high value of $\tau$ before stepping up.

   {
Despite the good performances shown in the previous results, we should not forget that the branch-and-cut procedure implemented in CPLEX is an exact algorithm, and therefore its usage as a heuristic is limited. To this aim, we report in Tables \ref{tab:cplex_longer_time_comparison} and \ref{tab:cplex_longer_time_comparison_poisson} the objective values obtained by CPLEX giving more time on the instances of the Gaussian and the Poissone set respectively. As the time limit increases the values obtained by the exact solver become increasingly lower. Of course, the larger the instance, the larger the branching tree, which implies a longer time for the branch-and-cut procedure. As it is clear from the table, the decay of our solution generalizes better in terms of the computational time needed to reach that value, meaning our procedure could be implemented as a starting point method to enhance the performances of the branch-and-cut procedure implemented in CPLEX.}

\begin{table}[H]
\centering
\footnotesize
			 \scalebox{0.92}{
			 \renewcommand\arraystretch{1.5} 
       {
\begin{tabular}{ccccccccccc}
\hline
J $\times$ M   & \#instances & \multicolumn{2}{c}{Mean}          & \multicolumn{2}{c}{Std dev} & \multicolumn{2}{c}{Max} & \multicolumn{2}{c}{Min} & Avg $\rho$ $(\%)$ \\ \hline
               &             & RL              & BC-CPLEX        & RL          & BC-CPLEX      & RL        & BC-CPLEX    & RL        & BC-CPLEX    &                   \\ \cline{3-10}
30 $\times$ 25 & 100         & 4785.3          & \textbf{3722.3} & 111.7       & 82.3          & 5107.3    & 3915.0      & 4579.7    & 3502.5      & -28.6             \\
35 $\times$ 30 & 100         & 5462.6          & \textbf{4729.0} & 115.1       & 92.7          & 5702.0    & 4906.9      & 5180.4    & 4468.4      & -15.6             \\
40 $\times$ 35 & 100         & 5989.5          & \textbf{5723.9} & 134.8       & 143.0         & 6349.8    & 6367.7      & 5715.4    & 5410.8      & -4.7              \\
45 $\times$ 40 & 100         & \textbf{6931.4} & 10429.2         & 130.2       & 17501.5       & 7304.9    & 89850.7     & 6599.8    & 6007.8      & -3.2              \\
50 $\times$ 45 & 100         & \textbf{8663.5} & 47061.5         & 147.4       & 51146.6       & 9015.6    & 125070.6    & 8198.9    & 7286.8      & \textbf{28.3}     \\ 
\end{tabular}}}
    \smallskip
	\caption{{RL vs the branch-and-cut algorithm
implemented in the known commercial solver CPLEX} on the {Gaussian set}: comparison on the objective value when the time limit is set to 60 seconds.}
	\label{tab:cplex_longer_time_comparison}
\end{table}

\begin{table}[H]
\centering
\footnotesize
			 \scalebox{0.92}{
			 \renewcommand\arraystretch{1.5} 
       {
\begin{tabular}{ccccccccccc}
\hline
J $\times$ M   & \#instances & \multicolumn{2}{c}{Mean}          & \multicolumn{2}{c}{Std dev} & \multicolumn{2}{c}{Max} & \multicolumn{2}{c}{Min} & Avg $\rho$ $(\%)$ \\ \hline
               &             & RL              & BC-CPLEX        & RL          & BC-CPLEX      & RL        & BC-CPLEX    & RL        & BC-CPLEX    &                   \\ \cline{3-10}
30 $\times$ 25 & 100         & 4167.1          & \textbf{3602.3} & 93.7        & 72.4          & 4404.0    & 3784.0      & 3939.0    & 3395.0      & -15.7             \\
35 $\times$ 30 & 100         & 5249.4          & \textbf{4655.6} & 136.0       & 63.0          & 5566.0    & 4814.0      & 4963.0    & 4488.0      & -12.8             \\
40 $\times$ 35 & 100         & 5947.5          & \textbf{5637.0} & 111.9       & 99.6          & 6194.0    & 5852.0      & 5735.0    & 5405.0      & -5.5              \\
45 $\times$ 40 & 100         & \textbf{7815.5} & 10683.9         & 182.6       & 17056.9       & 8263.0    & 87381.0     & 7453.0    & 6283.0      & -10.2             \\
50 $\times$ 45 & 100         & \textbf{8032.1} & 59159.0         & 181.9       & 54790.0       & 8410.0    & 128595.0    & 7697.0    & 7007.0      & \textbf{41.1}     \\
\end{tabular}}}
    \smallskip
	\caption{{RL vs the branch-and-cut algorithm
implemented in the known commercial solver CPLEX} on the {Poisson set}: comparison on the objective value when the time limit is set to 60 seconds.}
	\label{tab:cplex_longer_time_comparison_poisson}
\end{table}

   {Finally, to complete the analysis we performed an analysis on the standard Taillard's Benchmark instances derived from \cite{taillard1993benchmarks}, comparing our approach with the branch-and-cut procedure implemented in CPLEX setting the time limit to one and then to five minutes. We report these results in Table \ref{tab:cplex_tai}. As expected, given the small size of the instances, the values obtained by the commercial solver after one and five minutes of computation are dominant to our approach. For completeness, the branch-and-cut procedure implemented in CPLEX terminated at the time limit for every instance, both when the limit was set to one and to five. }

\begin{table}[H]
\centering
\footnotesize
			 \renewcommand\arraystretch{1.5} 
   {
\begin{tabular}{ccccccccccccccc}
\hline
J $\times$ M       & $\#$ instances & Statistics & \multicolumn{3}{c}{Algorithm}                                              \\ \hline
                   &                &           & BB-CPLEX$^1$ & BB-CPLEX$^5$ & RL     \\ \cline{4-6} 
15 $\times$ 15 & 10           & Mean      & 1286.2                          & \textbf{1253.4}                 & 1488.9 \\
Ta01-10            &                & Std. Dev  & 41.7                            & 33.3                            & 49.0   \\
                   &                & Max       & 1359.0                          & 1313.0                          & 1590.0 \\
                   &                & Min       & 1190.0                          & 1181.0                          & 1439.0 \\\hline
20 $\times$ 15 & 10           & Mean      & 1553.7                          & \textbf{1479.0}                 & 1735.5 \\
Ta11-20            &                & Std. Dev  & 53.9                            & 64.0                            & 52.9   \\
                   &                & Max       & 1647.0                          & 1619.0                          & 1822.0 \\
                   &                & Min       & 1437.0                          & 1405.0                          & 1657.0 \\\hline
20 $\times$ 20 & 10           & Mean      & 1859.7                          & \textbf{1744.7}                 & 2094.6 \\
Ta21-30            &                & Std. Dev  & 55.1                            & 48.8                            & 114.4  \\
                   &                & Max       & 1960.0                          & 1851.0                          & 2312.0 \\
                   &                & Min       & 1789.0                          & 1686.0                          & 1968.0 \\\hline
30 $\times$ 15 & 10           & Mean      & 2318.5                          & \textbf{2194.0}                 & 2370.6 \\
Ta31-40            &                & Std. Dev  & 111.1                           & 87.8                            & 164.6  \\
                   &                & Max       & 2524.0                          & 2350.0                          & 2670.0 \\
                   &                & Min       & 2102.0                          & 1985.0                          & 2075.0 \\\hline
30 $\times$ 20 & 10           & Mean      & 2709.8                          & \textbf{2496.6}                 & 2764.1 \\
Ta41-50            &                & Std. Dev  & 114.7                           & 91.7                            & 147.9  \\
                   &                & Max       & 2949.0                          & 2620.0                          & 3080.0 \\
                   &                & Min       & 2526.0                          & 2347.0                          & 2511.0 \\ \hline
\end{tabular}}
    \smallskip
	\caption{RL vs BB-CPLEX on the {Taillard benchmark instances}: comparison on the  on the makespan $C_{\max}$.  BB-CPLEX$^1$ represent the branch-and-cut procedure implemented in CPLEX with a time limit of one minute, while BB-CPLEX$^5$ stops after five minutes. }
	\label{tab:cplex_tai}
\end{table}

\subsubsection{Comparison with the Adaptive Algorithm} \label{sec:adaptive}  



As anticipated in the previous sections, we conducted a further comparison with the adaptive algorithm (ADA) described in \cite{Sotskov1996Adaptive} (see also \cite{Sotskov2013Adaptive}). The approach uses the weighted mixed disjunctive graph model to represent job shop scheduling problems and the conflict resolution strategy is applied to build a feasible schedule. Trained on a sample of job shop problems, the aim of the method is to produce knowledge on a benchmark of priority dispatching rules in order to solve similar large-scale job-shop problems by applying, by analogy, the acquired knowledge. In more detail, the adaptive algorithm proposed is characterized by two phases: a learning and an examination stage. In the former, the training instances are solved by an exact or approximate algorithm. Accordingly, information on successful decisions on conflict situations, represented by disjunctive edges in the graph, is stored in a learning table. In particular, for each resolved conflict edge, the algorithm gathers some characteristics based on priority dispatching rules and, along with the orientation of the edge, it saves all the computed information in a row of the learning database. In this way, it is possible to extract from the table a composite decision rule to produce a comprehensive specific heuristic.
In the examination stage, the adaptive scheduler solves new unseen instances by adopting decisions based on the derived heuristic.

As regard the learning stage, we trained the algorithm on the same instances reported in Table \ref{tab:training_sets} and we solved them using the branch-and-cut algorithm implemented in CPLEX, by setting a time limit equal to 420 seconds per instance. The time limit is a physical requirement since the training instances were many and the total training process would have been too time-consuming. During the learning, we used as benchmark priority rules the SPT (\textit{Shortest Processing Time}), LPT (\textit{Longest Processing Time}, FIFO (\textit{First In First Out}) and ECT (\textit{Earliest Completion Time}) heuristics.
In the examination stage, we tested the algorithm on 
unseen instances, described in Table \ref{tab:test_sets}, and, for each instance, 
due to the large amount of time required for the training, we randomly sampled 100 rows of the learning table.


In Tables \ref{tab:res_gauss_ada_obj}, \ref{tab:res_ada_poiss_obj}, \ref{tab:res_gauss_ada_time} and \ref{tab:res_poiss_ada_time}, we compared the adaptive algorithm with our approach on the Gaussian and Poisson test sets
We provide $\rho$ and $\tau$ values representing the percentage of the average improvement in terms of makespan and computational time respectively, where the single improvements of RL upon ADA on the $i$-th instance are computed as follows:
\begin{itemize}
    \item $\rho_i = \frac{C_{\max_{ADA_i}}-C_{\max_{RL_i}}}{{C_{\max_{ADA_i}}}}$
    \item $\tau_i = \frac{time_{ADA_i}-time_{RL_i}}{time_{ADA_i}}$
\end{itemize}
As we can see, our approach outperforms the adaptive algorithm on the test sets and the performances are very similar between the two set of instances. On the Gaussian set, the $\rho$ value is between 3\% and  27.5$\%$ and increases as the problems become more complex. 
Similarly, in the Poisson case, it ranges from 15$\%$ up to 26.3$\%$.
Concerning computational times, our RL approach outperforms the adaptive algorithm presenting an improvement in a range between 66$\%$ and 86$\%$ on both the test sets.

Finally, in Figure \ref{fig:performance_profiles_ADA}, we report the performance profiles for the Gaussian and the Poisson set, considering the computational time analysis first, and the objective value afterwords.
In particular, Figure \ref{fig:gaussianquality_ADA} and Figure \ref{fig:poissonquality_ADA} show that our approach beats more than 90\% of all the makespan returned by the ADA (both for gaussian and for poisson distribution). For Figure \ref{fig:gaussiantime_ADA} and \ref{fig:poissontime_ADA}, the situation is even more emphasized, with our approach faster than almost 100\% of all the times produced by the ADA, independently from the distribution.

\begin{table}[H]
\centering
\footnotesize
			 \renewcommand\arraystretch{1.5} 
\begin{tabular}{ccccccccccc}
\hline
J $\times$ M   & $\#$ instances & \multicolumn{2}{c}{Mean}                   & \multicolumn{2}{c}{Std. Dev} & \multicolumn{2}{c}{Max} & \multicolumn{2}{c}{Min} & Avg $\rho$ (\%)                \\ \hline
               &                & RL                               & ADA     & RL           & ADA         & RL         & ADA        & RL         & ADA        &                                \\\cline{3-10}
30 $\times$ 25 & 100            & \textbf{4785.3} & 4961.1  & 111.7        & 364.0       & 5107.3     & 6166.9     & 4579.7     & 4263.8     & \textbf{3.0}  \\
35 $\times$ 30 & 100            & \textbf{5462.6} & 6273.9  & 115.1        & 500.0       & 5702.0     & 8063.0     & 5180.4     & 5127.1     & \textbf{12.4} \\
40 $\times$ 35 & 100            & \textbf{5989.5} & 7698.5  & 134.8        & 535.0       & 6349.8     & 9312.9     & 5715.4     & 6277.4     & \textbf{21.8} \\
45 $\times$ 40 & 100            & \textbf{6931.4} & 9414.1  & 130.2        & 596.2       & 7304.9     & 11075.4    & 6599.8     & 8262.2     & \textbf{26.1} \\
50 $\times$ 45 & 100            & \textbf{8663.5} & 12007.4 & 147.4        & 883.6       & 9015.6     & 14312.1    & 8198.9     & 10351.4    & \textbf{27.5}
\end{tabular}
    \smallskip
	\caption{{RL vs ADA on the {Gaussian set}}: comparison on the makespan $C_{\max}$.}
	\label{tab:res_gauss_ada_obj}
\end{table}

\begin{table}[H]
\centering
\footnotesize
			 \renewcommand\arraystretch{1.5} 
\begin{tabular}{ccccccccccc}
\hline
J $\times$ M   & $\#$ instances & \multicolumn{2}{c}{Mean}  & \multicolumn{2}{c}{Std. Dev} & \multicolumn{2}{c}{Max} & \multicolumn{2}{c}{Min} & Avg $\rho$ (\%) \\ \hline
               &                & RL              & ADA     & RL           & ADA         & RL         & ADA        & RL         & ADA        &                 \\ \cline{3-10}
30 $\times$ 25 & 100            & \textbf{4167.1} & 4922.3  & 93.7         & 300.5       & 4404.0     & 5718.0     & 3939.0     & 4108.0     & \textbf{15.0}   \\
35 $\times$ 30 & 100            & \textbf{5249.4} & 6115.7  & 136.0        & 368.3       & 5566.0     & 7098.0     & 4963.0     & 5347.0     & \textbf{13.8}   \\
40 $\times$ 35 & 100            & \textbf{5947.5} & 7646.7  & 111.9        & 677.5       & 6194.0     & 9977.0     & 5735.0     & 6607.0     & \textbf{21.7}   \\
45 $\times$ 40 & 100            & \textbf{7815.5} & 9653.4  & 182.6        & 745.9       & 8263.0     & 11898.0    & 7453.0     & 7809.0     & \textbf{18.6}   \\
50 $\times$ 45 & 100            & \textbf{8032.1} & 10954.7 & 181.9        & 780.4       & 8410.0     & 13174.0    & 7697.0     & 9442.0     & \textbf{26.3}  
\end{tabular}
    \smallskip
	\caption{{RL vs ADA on the {Poisson set}}: comparison on the makespan $C_{\max}$.}
	\label{tab:res_ada_poiss_obj}
\end{table}

\begin{table}[H]
\centering
\footnotesize
			 \renewcommand\arraystretch{1.5} 
\begin{tabular}{ccccccccccc}
\hline
J $\times$ M   & $\#$ instances & \multicolumn{2}{c}{Mean} & \multicolumn{2}{c}{Std. Dev} & \multicolumn{2}{c}{Max} & \multicolumn{2}{c}{Min} & Avg $\tau$ (\%) \\ \hline
               &                & RL              & ADA    & RL           & ADA         & RL         & ADA        & RL         & ADA        &                 \\ \cline{3-10}
30 $\times$ 25 & 100            & \textbf{2.5}    & 7.3    & 0.3          & 0.5         & 3.5        & 9.2        & 1.9        & 6.7        & \textbf{66.2}   \\
35 $\times$ 30 & 100            & \textbf{4.6}    & 17.1   & 0.5          & 0.5         & 5.6        & 18.7       & 3.5        & 16.5       & \textbf{73.2}   \\
40 $\times$ 35 & 100            & \textbf{7.0}    & 29.7   & 1.0          & 0.7         & 10.3        & 32.0       & 5.7        & 28.5       & \textbf{76.3}   \\
45 $\times$ 40 & 100            & \textbf{11.8}   & 72.7   & 1.2          & 4.7         & 14.0       & 100.3      & 9.8        & 67.2       & \textbf{83.8}   \\
50 $\times$ 45 & 100            & \textbf{18.7}   & 138.4  & 1.9          & 2.3         & 24.4       & 147.8      & 15.9       & 129.1      & \textbf{86.5}  
\end{tabular}
    \smallskip
	\caption{{RL vs ADA on the {Gaussian set}}: comparison on the computational time (seconds).}
	\label{tab:res_gauss_ada_time}
\end{table}

\begin{table}[H]
\centering
\footnotesize
			 \renewcommand\arraystretch{1.5} 
\begin{tabular}{ccccccccccc}
\hline
J $\times$ M   & $\#$ instances & \multicolumn{2}{c}{Mean} & \multicolumn{2}{c}{Std. Dev} & \multicolumn{2}{c}{Max} & \multicolumn{2}{c}{Min} & Avg $\tau$ (\%) \\ \hline
               &                & RL              & ADA    & RL           & ADA         & RL         & ADA        & RL         & ADA        &                 \\ \cline{3-10}
30 $\times$ 25 & 100            & \textbf{2.4}    & 7.5    & 0.4          & 0.4         & 4.7        & 8.7        & 1.8        & 6.5        & \textbf{68.4}   \\
35 $\times$ 30 & 100            & \textbf{4.4}    & 16.5   & 0.5          & 0.7         & 6.1       & 18.8       & 3.4        & 15.6       & \textbf{73.4}   \\
40 $\times$ 35 & 100            & \textbf{7.2}    & 31.1   & 0.9          & 0.9         & 10.7        & 33.9       & 6.0        & 29.1       & \textbf{76.9}   \\
45 $\times$ 40 & 100            & \textbf{11.3}   & 68.3   & 1.1          & 2.0         & 14.1       & 75.0       & 9.7        & 63.6       & \textbf{83.4}   \\
50 $\times$ 45 & 100            & \textbf{17.1}   & 121.9  & 1.2          & 3.0         & 20.2       & 128.3      & 15.0       & 115.4      & \textbf{86.0}  
\end{tabular}
    \smallskip
	\caption{{RL vs ADA on the {Poisson set}}: comparison on the computational time (seconds).}
	\label{tab:res_poiss_ada_time}
\end{table}

\begin{figure}
\centering 
\footnotesize
\begin{subfigure}{0.5\textwidth}
\centering 
  \input{gaussian_time_ADAPTIVE}
    \caption{Gaussian set, computational time.}
    \label{fig:gaussiantime_ADA}
\end{subfigure}%
\hfill
\begin{subfigure}{0.5\textwidth}
  \centering 
  \input{poisson_time_ADAPTIVE}
    \caption{Poisson set, computational time.}
    \label{fig:poissontime_ADA}
\end{subfigure}

\vspace{.7cm}

\begin{subfigure}{0.5\textwidth}
\centering 
  \input{gaussian_obj_value_ADAPTIVE}
    \caption{Gaussian set, objective value.}
    \label{fig:gaussianquality_ADA}
\end{subfigure}%
\hfill
\begin{subfigure}{0.5\textwidth}
  \centering 
  \input{poisson_obj_value_ADAPTIVE}
    \caption{Poisson set, objective value.}
    \label{fig:poissonquality_ADA}
\end{subfigure}

\caption{Performance profiles.}
\label{fig:performance_profiles_ADA}
\end{figure}

\subsubsection{Comparison with scheduling heuristics based on Priority Rules} \label{sec:priority_rules} 

We compared our algorithm with 17 well-known heuristic rules for the JSSP, as described in \cite{Panwalkar1977ASurvey}.
These heuristics are based on a priority rule that selects the job to be processed.
In Table \ref{tab:heuristics_description}, we provide a list of the addressed heuristics with a short description of the rule. 

\begin{table}[H]
\centering
\footnotesize
			 \scalebox{0.92}{
			 \renewcommand\arraystretch{1.5} 
\begin{tabular}{ll}
\hline
\textbf{Priority rule} & \multicolumn{1}{c}{Description}                                                                                                          \\ \hline
SPT           & Select the job with the shortest processing time                                                                      \\
LPT           & Select the job with the longest processing time                                                                       \\
SSO           & Select the job with the shortest processing time of subsequent operation                                              \\
LSO           & Select the job with the longest processing time of subsequent operation                                               \\
SRM           & Select the job with the shortest remaining processing time not including the processing time of the current operation \\
LRM           & Select the job with the longest remaining processing time not including the processing time of the current operation  \\
FOPNR         & Select the job with fewest remaining operations                                                                       \\
SPT$+$SSO     & Select the job with the minimum sum of the processing times of the current and subsequent operation                    \\
LPT$+$LSO     & Select the job with the maximum sum of the processing times of the current and subsequent operation                    \\
SPT*TWK       & Select the job with the minimum product of current operation processing time and total working time                   \\
LPT*TWK       & Select the job with the maximum product of current operation processing time and total working time                   \\
SPT/TWK       & Select the job with the minimum ratio of current operation processing time to total working time                      \\
LPT/TWK       & Select the job with the maximum ratio of current operation processing time to total working time                      \\
SPT*TWKR      & Select the job with the minimum product of current operation processing time to total remaining working time          \\
LPT*TWKR      & Select the job with the maximum product of current operation processing time and total remaining working time         \\
SPT/TWKR      & Select the job with the minimum ratio of current operation processing time to total remaining working time            \\
LPT/TWKR      & Select the job with the minimum ratio of current operation processing time to total remaining working time            \\ 
\end{tabular}}
    \smallskip
	\caption{PDR-based heuristics}
	\label{tab:heuristics_description}
\end{table}


We evaluated our method against the 17 PR-based heuristics on the {Gaussian} and {Poisson} test sets already described in the previous sections. The results are reported in tables \ref{tab:res_gauss_obj_heuristic} and \ref{tab:res_poiss_obj_heuristic}, respectively. 
For each problem size, we report the average makespan $C_{\max}$ over all the $100$ instances of the same size.


The only heuristic that performs better than our RL approach is the LRM (\textit{Longest Remaining Machining time}). More in detail, the LRM heuristic defeats RL on the $30\times25 $ and $50\times 45$ instances of the {Gaussian set}, reaching an improvement of around 1.7\% and 5.7\% respectively. Concerning the {Poisson set}, LRM overcomes RL only on the $45 \times 40$ instances, with an improvement of around 2\%.
The rules LPT*TWKR, SPT/TWKR and LPT/TWKR return average makespans around 2 times higher than our approach, while all the other heuristics produces makespans 5 times or higher. 
Therefore we may say that the proposed RL method has dominant performances compared to the selected heuristics in almost every class of test instances.


\begin{landscape}
\thispagestyle{empty} 

\begin{table}
\centering
\footnotesize
			 \scalebox{0.7}{
			 \renewcommand\arraystretch{1.5} 
\begin{tabular}{cccccccccccccccccccc}
\hline
J $\times$ M   & $\#$ instances & \multicolumn{18}{c}{\textbf{Mean makespan $\hat C_{\max}$}}                                                                                                                                                  \\ \hline
               &             & RL             & SPT     & LPT     & SSO     & LSO     & SRM      & LRM             & FOPNR    & SPT+SSO & LPT+LSO & SPT*TWK & LPT*TWK & SPT/TWK & LPT/TWK & SPT*TWKR & LPT*TWKR & SPT/TWKR & LPT/TWKR \\
\cline{3-20}
30 $\times$ 25 & 100         & 4785.2          & 27043.3 & 27522.2 & 27239.2 & 27288.8 & 39828.6  & \textbf{4704.8} & 40818.6  & 27207.9 & 25069.4 & 28159.7 & 28261.6 & 28015.4 & 28062.6 & 38418.9  & 7051.0   & 6461.5   & 38121.0  \\
35 $\times$ 30 & 100         & \textbf{5462.5} & 38957.3 & 39576.3 & 39449.8 & 39629.7 & 58356.8  & 5613.1          & 57553.7  & 39286.9 & 36663.2 & 40637.9 & 39968.7 & 39989.2 & 40638.7 & 55808.8  & 9586.9   & 8581.0   & 55593.6  \\
40 $\times$ 35 & 100         & \textbf{5989.4} & 50882.5 & 51497.9 & 50997.7 & 51235.2 & 74920.2  & 6417.3          & 75766.9  & 50176.3 & 48228.3 & 51956.8 & 52535.2 & 51826.4 & 52398.5 & 72103.3  & 11405.3  & 10155.9  & 71926.3  \\
45 $\times$ 40 & 100         & \textbf{6931.4} & 67592.2 & 68130.0 & 67756.5 & 67755.1 & 101495.4 & 7181.2          & 97722.8  & 67400.7 & 64684.4 & 69874.5 & 69997.1 & 69245.4 & 69943.0 & 96495.7  & 14598.4  & 12757.1  & 96382.9  \\
50 $\times$ 45 & 100         & 8663.5          & 88020.5 & 87741.0 & 88539.9 & 87247.4 & 126439.2 & \textbf{8168.0} & 127350.1 & 87411.1 & 84192.0 & 89840.5 & 89369.3 & 89151.4 & 89384.8 & 122160.7 & 18433.8  & 15937.8  & 122109.3 \\ 
\end{tabular}}
    \smallskip
	\caption{{RL vs HEURISTICS on the {Gaussian set}}: results for the objective value, i.e. makespan $C_{\max}$.}
	\label{tab:res_gauss_obj_heuristic}
\end{table}


\begin{table}
\centering
\footnotesize
			 \scalebox{0.7}{
			 \renewcommand\arraystretch{1.5} 
\begin{tabular}{cccccccccccccccccccc}
\hline
J $\times$ M   & \#instances & \multicolumn{18}{c}{\textbf{Mean Makespan $\hat C_{\max}$}}                                                                                                                                                  \\ \hline
               &             & RL             & SPT     & LPT     & SSO     & LSO     & SRM      & LRM             & FOPNR    & SPT+SSO & LPT+LSO & SPT*TWK & LPT*TWK & SPT/TWK & LPT/TWK & SPT*TWKR & LPT*TWKR & SPT/TWKR & LPT/TWKR \\
               \cline{3-20}
30 $\times$ 25 & 100         & \textbf{4167.1} & 26094.0                 & 27272.4                 & 26676.2                 & 27139.1                 & 39473.4                 & 4504.8                                 & 38448.4                   & 26498.9                     & 24493.9                     & 27263.5                     & 27263.5                     & 27120.9                     & 27445.1                     & 37830.9                      & 6749.6                       & 6247.8                       & 37762.4                      \\
35 $\times$ 30 & 100         &\textbf{5249.3} & 37867.5                 & 38509.7                 & 38401.7                 & 38075.3                 & 56306.4                 & 5379.0                                 & 56190.3                   & 37894.1                     & 35910.0                     & 39368.3                     & 39368.3                     & 38904.1                     & 39545.2                     & 53946.5                      & 9092.5                       & 8290.9                       & 54070.7                      \\
40 $\times$ 35 & 100         &\textbf{5947.4} & 52884.2                 & 53281.0                 & 53186.9                 & 53454.9                 & 77287.6                 & 6225.5                                 & 76786.0                   & 52041.8                     & 50473.1                     & 53645.9                     & 53645.9                     & 53352.1                     & 54144.3                     & 74676.3                      & 11563.9                      & 10244.0                      & 74043.8                      \\
45 $\times$ 40 & 100         &  7815.5          & 67844.4                 & 68982.8                 & 67885.8                 & 69024.4                 & 99289.5                 &  \textbf{7655.5} & 98925.0                   & 68173.4                     & 65585.7                     & 69691.3                     & 69691.3                     & 69601.5                     & 70168.4                     & 94933.3                      & 14873.2                      & 13379.1                      & 94795.6                      \\
50 $\times$ 45 & 100         &  \textbf{8032.1} & 85492.6                 & 87243.2                 & 85696.6                 & 86781.4                 & 125574.0                & 8085.7          & 123372.0                  & 85857.1                     & 83168.9                     & 87738.4                     & 87738.4                     & 87246.0                     & 88250.1                     & 120762.4                     & 17430.5                      & 15500.7                      & 120368.1                     \\
\end{tabular}}
    \smallskip
	\caption{{RL vs HEURISTICS on the {Poisson set}}: results for the objective value, i.e. makespan $C_{\max}$.}
	\label{tab:res_poiss_obj_heuristic}
\end{table}
\end{landscape}

\subsubsection{Comparison with the Reinforcement Learning algorithm in \cite{zhang2020learning} } \label{sec:zhang}
   {
Finally, we compared our approach with the one presented in \cite{zhang2020learning}. In this work,  the authors proposed a method based on Graph Neural Networks, a special class of Deep Neural Network models dealing with graph-structured data (see \cite{wu2020comprehensive} and \cite{zhou2020graph} for a comprehensive review). From the result in \cite{blazewicz2000disjunctive}, the JSSP can be formulated as a disjunctive graph and fed to the network. Differently from us, they used the PPO algorithm using a clipping loss to compute the policy. }

   {
The training and testing phases were performed by generating instances with the Taillard method \cite{taillard1993benchmarks}. The authors created different classes of problems: $6 \times 6$, $10 \times 10$, $15 \times 15$, $20 \times 15$, $20 \times 20$, and $30 \times 20$. Each class had $100$ instances. The authors trained a different model for each class, performing $10000$ PPO iterations using $4$ trajectories each.
}

   {To compare our approach to the one in \cite{zhang2020learning}, we re-trained our model on the same Taillard instances, excluding the classes  $30 \times 15$ and  $30 \times 20$ for computational bounds on our resources. Since the training of our model can handle different JSSP instances with heterogeneous sizes at the same time, we preferred to perform a one-shot training using the sets of instances ranging from $6 \times 6$ to $20 \times 20$. We performed $25000$ PPO iterations with $5$ trajectories each.}

   {
We reported our results in Tables \ref{tab:zhang} and \ref{tab:zhang2}. 
In the first Table, we show the value obtained in terms of objective on the Taillard's benchmark set of known Taillard instances taken from \cite{taillard1993benchmarks} (see Table \ref{tab:taillard}).
The results show that our model is able to slightly beat the approach proposed in \cite{zhang2020learning}. This is surprising since the nature of our training tends to be less specialized on the single group of instances, as opposed to \cite{zhang2020learning} where the results are generated using models trained only on the class at hand. The classes $30 \times 15$ and $30 \times 20$ are included to show the quality of our solution which was not trained on those distributions, differently from the results proposed in \cite{zhou2020graph}, which have specialized models for the $30 \times 20$ case. 
}

\begin{table}[H]
\centering
\footnotesize
			\scalebox{1.10}{
 \renewcommand\arraystretch{1.5}
   {\begin{tabular}{cccccccc}\hline
J $\times$ M   & Instance IDs & $\#$ instances & \multicolumn{4}{c}{Our RL}                     & Reported on \cite{zhang2020learning} \\ \hline
               &         &  & Mean            & Std. Dev & Max    & Min    & Mean                                            \\ \cline{4-8}
15 $\times$ 15 & Ta01-10 & 10  & \textbf{1488.9} & 49.0     & 1590.0 & 1439.0 & 1547.7                                          \\
20 $\times$ 15 & Ta11-20  & 10  & \textbf{1735.5} & 52.9     & 1822.0 & 1657.0 & 1774.7                                          \\
20 $\times$ 20 & Ta21-30 & 10   & \textbf{2094.6} & 114.4    & 2312.0 & 1968.0 & 2128.1                                          \\
30 $\times$ 15 & Ta31-40 & 10   & \textbf{2370.6} & 164.6    & 2670.0 & 2075.0 & 2378.8                                          \\
30 $\times$ 20 & Ta41-50 & 10   & 2764.1          & 147.9    & 3080.0 & 2511.0 & \textbf{2603.9}                                \\
\end{tabular}}}
    \smallskip
	\caption{RL vs \cite{zhang2020learning} on the Taillard benchmark instances: comparison on the makespan $C_{\max}$.}
	\label{tab:zhang}
\end{table}

   { In the second Table (\ref{tab:zhang2}), we compared the two approaches on newly generated instances, obtained with the Taillard method \cite{taillard1993benchmarks}, using the generator provided in \cite{zhang2020learning}. In the Table, we re-run the different models available from \cite{zhang2020learning}, and we computed the objective values for every new class. In particular, we created a total of $6$ classes with $100$ instances each having the following specifics:  $6 \times 6$, $10 \times 10$, $15 \times 10$, $15 \times 15$, $20 \times 10$, and $20 \times 20$. Our approach is capable of beating different variations of the algorithm presented in \cite{zhang1995reinforcement} up to the size of $30 \times 20$, except for the instances in the class $20 \times 20$. }

   {We included an analysis on timing in the Appendix, showing the computational times are comparable, with no relevant difference between the two approaches.}

\begin{table}[H]
\centering
\footnotesize
			 \renewcommand\arraystretch{1.5} 
   {\begin{tabular}{ccccccccc}\hline
J $\times$ M   & $\#$ instances & Statistics & \multicolumn{6}{c}{Algorithm}                                                                                                                                                                                                   \\ \hline
               &                &           & \cite{zhang2020learning} $6 \times 6$ & \cite{zhang2020learning} $10 \times 10$ & \cite{zhang2020learning} $15 \times 15$ & \cite{zhang2020learning} $20 \times 20$ & \cite{zhang2020learning} $30 \times 20$ & Our RL          \\ \cline{4-9} 
$6 \times 6$   & 100            & Mean      & 581.9   & 574.5     & 571.7     & 573.2     & 570.8     & \textbf{544.7}  \\
               &                & Std. Dev  & 81.0    & 71.8      & 79.0      & 70.2      & 71.8      & 68.0            \\
               &                & Max       & 814.0   & 804.0     & 804.0     & 775.0     & 804.0     & 814.0           \\
               &                & Min       & 427.0   & 415.0     & 424.0     & 405.0     & 405.0     & 403.0           \\ \hline
$10 \times 10$ & 100            & Mean      & 1051.5  & 995.9     & 996.8     & 998.0     & 991.8     & \textbf{938.0}  \\
               &                & Std. Dev  & 96.2    & 76.2      & 78.9      & 74.3      & 75.9      & 71.4            \\
               &                & Max       & 1327.0  & 1187.0    & 1218.0    & 1242.0    & 1187.0    & 1133.0          \\
               &                & Min       & 842.0   & 862.0     & 820.0     & 871.0     & 781.0     & 770.0           \\ \hline
$15 \times 10$ & 100            & Mean      & 1306.9  & 1225.2    & 1229.0    & 1222.4    & 1222.0    & \textbf{1182.9} \\
               &                & Std. Dev  & 106.9   & 100.8     & 92.6      & 92.1      & 95.2      & 96.4            \\
               &                & Max       & 1573.0  & 1460.0    & 1504.0    & 1447.0    & 1440.0    & 1412.0          \\
               &                & Min       & 1064.0  & 995.0     & 1007.0    & 1012.0    & 995.0     & 968.0           \\ \hline
$15 \times 15$ & 100            & Mean      & 1636.6  & {1502.9}  & 1505.1    & 1503.3    & 1503.1    & \textbf{1441.8} \\
               &                & Std. Dev  & 103.8   & 105.1     & 107.1     & 96.3      & 97.6      & 98.9            \\
               &                & Max       & 1864.0  & 1841.0    & 1898.0    & 1785.0    & 1861.0    & 1752.0          \\
               &                & Min       & 1425.0  & 1321.0    & 1322.0    & 1311.0    & 1267.0    & 1261.0          \\ \hline
$20 \times 10$ & 100            & Mean      & 1567.2  & 1474.6    & {1470.4}  & 1476.2    & 1478.3    & \textbf{1436.4} \\
               &                & Std. Dev  & 104.9   & 96.6      & 95.8      & 98.5      & 104.8     & 108.8           \\
               &                & Max       & 1893.0  & 1682.0    & 1680.0    & 1755.0    & 1711.0    & 1756.0          \\
               &                & Min       & 1367.0  & 1201.0    & 1147.0    & 1188.0    & 1178.0    & 1223.0          \\ \hline
$20 \times 20$ & 100            & Mean      & 2216.2  & 1997.4    & 1993.2    & \textbf{1984.3}                         & 1996.7    & 2012.6          \\
               &                & Std. Dev  & 127.3   & 106.3     & 95.8      & 98.6      & 114.8     & 93.2            \\
               &                & Max       & 2528.0  & 2344.0    & 2253.0    & 2204.0    & 2355.0    & 2241.0          \\
               &                & Min       & 1954.0  & 1778.0    & 1722.0    & 1790.0    & 1775.0    & 1746.0          \\ \hline
\end{tabular}}
    \smallskip
	\caption{RL vs \cite{zhang2020learning} on the {Taillard generated instances}: comparison on the makespan $C_{\max}$.}
	\label{tab:zhang2}
\end{table}

\section{Conclusions}
In this paper, we investigated how to solve the Job Shop Scheduling problem (JSSP) through reinforcement learning, aiming to make the learning agent flexible for tackling instances with a variable number of jobs, tasks, and machines.

We first formulated the JSSP as a Markov Decision Process, which was fundamental to inscribe the problem in an actor-critic scheme. The method adopted takes inspiration from the Proximal Policy Optimization, \cite{schulman2017proximal},  using a dynamic adaptation of the penalty term to facilitate exploitation over exploration, and vice-versa, depending on the situation.

In the second phase, we studied several classes of deep models that could fit the JSSP, eventually landing on a double incident LSTM framework, where each LSTM works as a projection into a fixed space. The actor ends with an action masking to control feasibility, combined with a soft-max function
to recreate a discrete probability distribution, aka the policy estimator. At its bottom, the critic has an encoder network, collapsing the embeddings of the second LSTM into a scalar, representing the state-value function estimator. 

Our algorithm can generalize to a certain extent to instances with larger sizes, and with different distributions, than the one used in the training phase. The approach shows a decisive improvement towards the deterministic mixed-integer branch-and-cut algorithm {implemented in known solver CPLEX, the adaptive algorithm implemented in \cite{Sotskov1996Adaptive}, 17 priority rule-based heuristics, and a Deep Reinforcement learning algorithm, especially in terms of the makespan value, finally proving it is possible to generate new efficient greedy heuristics just from learning-based methodologies.}

{\setstretch{1}
\bibliographystyle{apalike}
\bibliography{biblio}}

\section*{Appendix}

\subsection*{Settings tables.}

\begin{table}[ht!]
	\centering
	\footnotesize
			\renewcommand\arraystretch{1.5} 	\begin{tabular}{cc}
			\hline
Parameter & 
Setting\\
\hline

\textit{LSTM 1} hidden size & $110$\\

\textit{LSTM 2} hidden size & $110\cdot 2$\\

\textit{FFN} number of hidden layers & $3$\\

\textit{FFN} input size & $110\cdot20$\\

\textit{FFN} 1st hidden layer size & $110\cdot 10$\\

\textit{FFN} 2nd hidden layer size & $110\cdot 5$\\

\textit{FFN} 3rd   hidden layer size & $110$\\

\textit{FFN} output layer size & $1$\\
\hline 
\end{tabular}
    \smallskip
	\caption{Actor-critic network configurations.}
	\label{tab:model_conf}
\end{table}

\begin{table}[ht!]
	\centering
		\footnotesize
		\renewcommand\arraystretch{1.5} 	\begin{tabular}{cccc}
\hline
Parameter & Setting \\
\hline
Number of episodes & $5000$\\
Roll-outs per episode  & $10$\\
Random choice prob.ty & From $\epsilon = 0.2$ down to $\epsilon=0$\\
KL penalty coefficient $\beta$ & $15$ \\
Target KL divergence $\delta$ & $0.05$\\
\hline
Actor optimizer & ADAM\\

Actor learning rate & $10^{-4}$ \\

Actor optimization steps & $1$\\

Critic optimizer & ADAM\\

Critic learning rate & $10^{-4}$\\

Critic optimization steps & $3$\\

Mini-batch size & $N\cdot T$\\
\hline
\end{tabular}
    \smallskip
	\caption{Training settings.}
	\label{tab:hyper_setting}
\end{table}

\subsection*{JSSP Mixed-integer formulation}

We used the following mixed-integer linear program for the JSSP in the solver CPLEX. 

$$
\left\{
\begin{array}{rll}
\min_{{\boldsymbol{t}},{\boldsymbol{x}}, C_{\textit{max}}} &  C_{\textit{max}}\\
\text{s.t.} & t_{jh}-t_{jk}\geq p_{jk},&\forall\  (j,k),(j,h)\in \mathcal{O},  (j,k)\prec (j,h)\\
& t_{jk}-t_{ik}\geq p_{ik}-M x_{jik},&\forall\  (j,k),(i,k)\in \mathcal{O}, j<i\\
& t_{ik}-t_{jk}\geq p_{jk}-M(1-x_{jik}),&\forall\  (j,k),(i,k)\in \mathcal{O}, j<i\\
& C_{\textit{max}}\geq t_{jk}+p_{jk},&\forall\ (j,k)\in \mathcal{O}\\
& t_{jk}\geq 0,&\forall\ (j,k)\in \mathcal{O}\\
& x_{jik}\in \{0,1\},&\forall\ (j,k),(i,k)\in \mathcal{O}, j<i
\end{array}
\right.
$$

Where:
\begin{itemize}
\item $p_{jk}$ are the processing times.
\item $M=\displaystyle \left\lceil \sum_{(j,k) \in \mathcal{O}} p_{jk} \right\rceil$ is a Big-M value.
\item $C_{\textit{max}}$ is  a continuous variable indicating the makespan.
    \item $t_{jk}$ are {continuous} variables indicating the non-negative starting time of operation $(j,k)$.
    \item $x_{jik}$ are {binary} variables defined as
    $$
    x_{jik}=\left\{ \begin{array}{l}
		    1, \text{\quad if job } j \text{ precedes } \text{job }i  \text{ on machine }k\\
		    0, \text{\quad otherwise} 
		\end{array}\right.\;
	$$
\end{itemize}

\subsection*{Table of times for the comparison in \cite{zhang2020learning}}

\begin{table}[H]
\centering
\footnotesize
			 \renewcommand\arraystretch{1.5} 
   {\begin{tabular}{ccccccccc}
J $\times$ M   & $\#$ instances  & Statistic & \multicolumn{6}{c}{Algorithm}                           \\ \hline
               &     &           & \cite{zhang2020learning} $6 \times 6$ & \cite{zhang2020learning} $10 \times 10$ & \cite{zhang2020learning} $15 \times 15$ & \cite{zhang2020learning} $20 \times 20$ & \cite{zhang2020learning} $30 \times 20$ & Our RL        \\ \cline{4-9} 
$6 \times 6$   & 100 & Mean      & 0.12                                                   & 0.12                                                     & 0.12                                                     & 0.12                                                     & 0.12                                                     & \textbf{0.06} \\
               &     & Std. Dev  & 0.08   & 0.05     & 0.05     & 0.04     & 0.04     & \textbf{0.01}          \\
               &     & Max       & 0.88   & 0.88     & 0.88     & 0.88     & 0.88     & \textbf{0.13}          \\
               &     & Min       & 0.10   & 0.10     & 0.10     & 0.09     & 0.09     & \textbf{0.04}          \\ \hline
$10 \times 10$ & 100 & Mean      & 0.31   & 0.32     & 0.32     & 0.32     & 0.32     & \textbf{0.17} \\
               &     & Std. Dev  & 0.08   & 0.06     & 0.05     & 0.04     & 0.04     & \textbf{0.01}          \\
               &     & Max       & 1.02   & 1.02     & 1.02     & 1.02     & 1.02     & \textbf{0.19}          \\
               &     & Min       & 0.27   & 0.27     & 0.27     & 0.27     & 0.27     & \textbf{0.14}          \\ \hline
$15 \times 10$ & 100 & Mean      & 0.45   & 0.49     & 0.48     & 0.49     & 0.49     & \textbf{0.31} \\
               &     & Std. Dev  & 0.07   & 0.08     & 0.07     & 0.07     & 0.07     & \textbf{0.01}         \\
               &     & Max       & 1.01   & 1.01     & 1.01     & 1.01     & 1.01     & \textbf{0.35}          \\
               &     & Min       & 0.37   & 0.37     & 0.37     & 0.37     & 0.37     & \textbf{0.30}          \\ \hline
$15 \times 15$ & 100 & Mean      & \textbf{0.67}                                          & 0.75     & 0.77     & 0.78     & 0.79     & 0.72          \\
               &     & Std. Dev  & \textbf{0.08}   & 0.14     & 0.14     & 0.14     & 0.13     & 0.13          \\
               &     & Max       & 1.30   & 1.47     & 1.47     & 1.47     & 1.47     & \textbf{1.14}          \\
               &     & Min       & 0.58   & 0.58     & 0.58     & 0.58     & 0.58     & \textbf{0.50}          \\ \hline
$20 \times 10$ & 100 & Mean      & \textbf{0.57}                                          & 0.63     & 0.65     & 0.67     & 0.68     & 0.72          \\
               &     & Std. Dev  & 0.07   & 0.10     & 0.10     & 0.12     & 0.12     & \textbf{0.06}          \\
               &     & Max       & 1.23   & 1.23     & 1.23     & 1.26     & 1.26     &\textbf{0.98}          \\
               &     & Min       & \textbf{0.53}   & \textbf{0.53}     & \textbf{0.53}     & \textbf{0.53}     & \textbf{0.53}     & 0.64          \\ \hline
$20 \times 20$ & 100 & Mean      & \textbf{1.42}                                          & 1.52     & 1.54     & 1.56     & 1.58     & 1.94          \\
               &     & Std. Dev  & \textbf{0.25}   & 0.26     & 0.28     & 0.27     & 0.31     & 0.38          \\
               &     & Max       & \textbf{2.00}   & 2.35     & 2.78     & 2.78     & 3.58     & 3.52          \\
               &     & Min       & \textbf{1.11}   & \textbf{1.11}     & \textbf{1.11}     & \textbf{1.11}     & \textbf{1.11}     & 1.64          \\ \hline
\end{tabular}
}
    \smallskip
	\caption{RL vs \cite{zhang2020learning} on the {Taillard generated instances}: comparison on the computational time.}
	\label{tab:zhang2_time}
\end{table}

\end{document}